%% file: jnl-2021-jmmct-sldm.tex
\documentclass[journal]{./IEEEtran}
\ifCLASSINFOpdf
\else
\fi
\input{./packages.tex}
\input{./definitions.tex}

\begin{document}
%
\title{A Single-Layer Dual-Mesh Boundary Element Method for Multiscale Electromagnetic Modeling\\ of Penetrable Objects in Layered Media}
%
%
%

\author{Shashwat~Sharma,~\IEEEmembership{Graduate Student Member,~IEEE,}
        and~Piero~Triverio,~\IEEEmembership{Senior Member,~IEEE}
\thanks{S. Sharma and P. Triverio are with the Edward S. Rogers Sr. Department of Electrical \& Computer Engineering, University of Toronto, Toronto, ON, M5S 3G4 Canada, e-mails: shash.sharma@mail.utoronto.ca, piero.triverio@utoronto.ca.}
\thanks{This work was supported by the Natural Sciences and Engineering Research Council of Canada (Collaborative Research and Development Grants program), by Advanced Micro Devices, and by CMC Microsystems.}
\thanks{Manuscript received $\ldots$; revised $\ldots$.}}

%
%

\markboth{IEEE Journal on Multiscale and Multiphysics Computational Techniques}%
{Sharma \MakeLowercase{\textit{et al.}}: Bare Demo of IEEEtran.cls for IEEE Journals}
%



\maketitle

\begin{abstract}
A surface integral representation of Maxwell's equations allows the efficient electromagnetic (EM) modeling of three-dimensional structures with a two-dimensional discretization, via the boundary element method (BEM).
However, existing BEM formulations either lead to a poorly conditioned system matrix for multiscale problems, or are computationally expensive for objects embedded in layered substrates.
This article presents a new BEM formulation which leverages the surface equivalence principle and Buffa-Christiansen basis functions defined on a dual mesh, to obtain a well-conditioned system matrix suitable for multiscale EM modeling.
Unlike existing methods involving dual meshes, the proposed formulation avoids the double-layer potential operator for the surrounding medium, which may be a stratified substrate requiring the use of an advanced Green's function.
This feature greatly alleviates the computational expense associated with the use of Buffa-Christiansen functions.
Numerical examples drawn from several applications, including remote sensing, chip-level EM analysis, and metasurface modeling, demonstrate speed-ups ranging from~$3\times$ to~$7\times$ compared to state-of-the-art formulations.
\end{abstract}

\begin{IEEEkeywords}
Maxwell's equations, boundary element method, integral equations, multiscale modeling, layered substrate.
\end{IEEEkeywords}

%
\IEEEpeerreviewmaketitle


\input{introduction.tex}

\input{methods.tex}

\input{results.tex}

\section{Conclusion}\label{sec:conclusion}
The electromagnetic modeling of complex structures with the boundary element method (BEM) requires efficient formulations, which can robustly handle multiple scales of feature size, operating frequency, and material properties of the objects.
In order to obtain a well-conditioned system matrix, existing formulations require the use of expensive Buffa-Christiansen basis functions and the double-layer potential operator associated with the background medium, which is often a layered substrate.
These requirements lead to a high computational cost for large problems.
This article introduces a novel single-layer formulation which leads to a well-conditioned system matrix for challenging multiscale problems, while avoiding the aforementioned double-layer operator.
The proposed method leverages the advantages of the Buffa-Christiansen functions while greatly alleviating the associated computational cost.
Consequently, speed-ups of~$3\times$ to~$7\times$ compared to existing methods are achieved for realistic structures, as demonstrated via several numerical examples.


%

%

\section*{Acknowledgment}

The authors would like to thank the anonymous reviewers for their thoughtful and constructive feedback.

\ifCLASSOPTIONcaptionsoff
  \newpage
\fi



\bibliographystyle{./IEEEtran}
\bibliography{./IEEEabrv,./bibliography}

\begin{IEEEbiography}[{\includegraphics[width=1in,height=1.25in,clip,keepaspectratio]{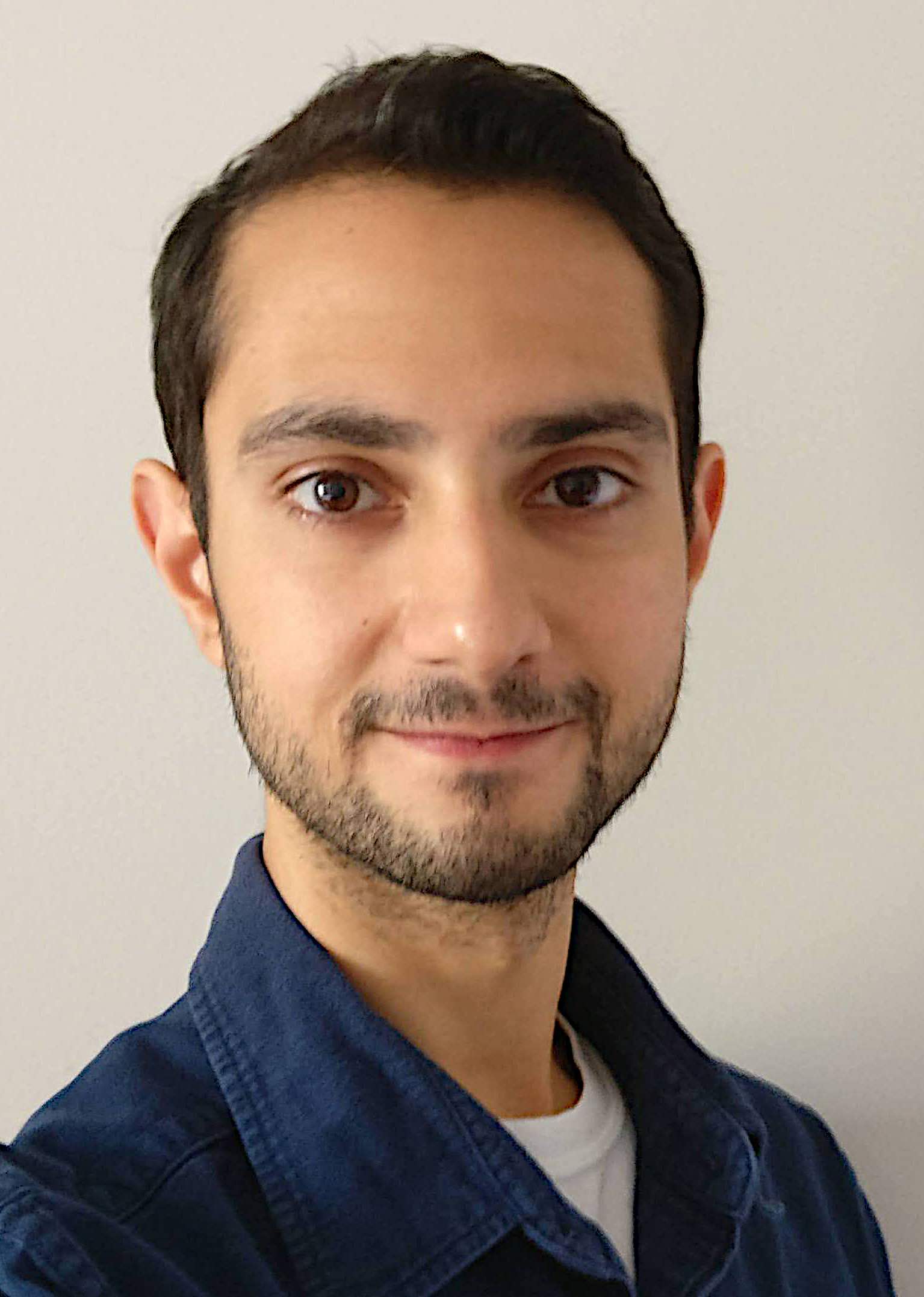}}]{Shashwat Sharma} (S'18)
	received the B.A.Sc. degree in engineering physics and the M.A.Sc. degree in electrical engineering from the University of Toronto, Canada, in 2014 and 2016, respectively. From September 2016 to March 2017 he worked as a computational science intern at Autodesk Research, Toronto. Since 2017 he has been working towards the Ph.D. degree in electrical engineering at the University of Toronto. His research focuses on computational electromagnetics, with an emphasis on fast and robust integral equation techniques for multiscale electromagnetic modeling. His interests include all aspects of numerical methods, computational science and high-performance computing in the context of electromagnetics.
	
	Mr.~Sharma received the Piergiorgio Uslenghi
	Letters Prize Paper Award (2021) and placed second at the CNC/USNC-URSI Student Paper Competition of the 2020 IEEE International Symposium on Antennas and Propagation and North American Radio Science Meeting (AP-S/URSI). He was also a finalist for the Best Student Paper Award at the 2018 IEEE Conference on Electrical
	Performance of Electronic Packaging and Systems. Mr.~Sharma has received honorable mentions for his contributions to the AP-S/URSI symposia in 2019 and 2020, and a TICRA-EurAAP Grant at the 2021 European Conference on Antennas and Propagation. He is currently serving as the vice chair for the University of Toronto student chapter of the IEEE Antennas and Propagation Society.
\end{IEEEbiography}

\begin{IEEEbiography}[{\includegraphics[width=1in,height=1.25in,clip,keepaspectratio]{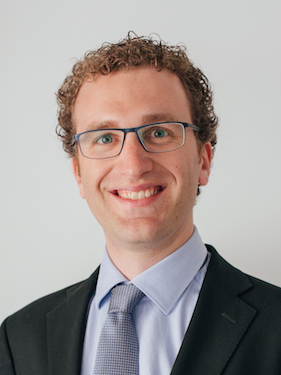}}]{Piero Triverio} (S'06 -- M'09 -- SM'16)
	received the Ph.D. degree in Electronic Engineering from Politecnico di
	Torino, Italy, in 2009. He is an Associate Professor in The Edward S.
	Rogers Sr. Department of Electrical \& Computer Engineering (ECE) at the
	University of Toronto, and in the Institute of Biomaterials and
	Biomedical Engineering (IBBME). He holds the Canada Research Chair in
	Computational Electromagnetics. His research interests include signal
	integrity, computational electromagnetism, model order reduction, and
	computational fluid dynamics applied to cardiovascular diseases.
	
	Prof.~Triverio and his students received the Piergiorgio Uslenghi
	Letters Prize Paper Award (2021), the Best Paper Award of the IEEE
	Transactions on Advanced Packaging (2007), the EuMIC Young Engineer
	Prize (2010), and the Ontario Early Researcher Award (2016). From 2013
	to 2018, Triverio held the Canada Research Chair in Modeling of
	Electrical Interconnects. Triverio and his students were awarded the
	Best Paper Award of the IEEE Conference on Electrical Performance of
	Electronic Packaging and Systems (2008, 2017), and several Best Student
	Paper Awards at international symposia. He serves as an Associate Editor
	for the IEEE TRANSACTIONS ON COMPONENTS, PACKAGING AND MANUFACTURING
	TECHNOLOGY. He is a member of the Technical Program Committee of the
	IEEE Workshop on Signal and Power Integrity, and of the IEEE Conference
	on Electrical Performance of Electronic Packaging and Systems.
\end{IEEEbiography}

\end{document}

%% file: packages.tex
\usepackage{cite}
\usepackage{amsmath}
\usepackage[pdftex]{graphicx} 
\usepackage{subfig}
\usepackage{array} 
\usepackage{fixltx2e}

\usepackage{booktabs} 
\usepackage{amssymb}
\usepackage{mathtools}
\usepackage{cancel} 
\usepackage{bm} 
\usepackage{esint} 
\usepackage{amsthm}
\usepackage{siunitx}
\usepackage{xfrac}

\usepackage{tikz}
\usepackage{tikz-3dplot}
\usetikzlibrary{backgrounds,3d,shadings,shapes.misc,decorations.pathmorphing,shapes,calc,fadings,shadows.blur}


%% file: definitions.tex
\usepackage{amssymb,amsmath}

\usepackage{bm}
\newcommand{\bs}[1]{\bm{\mathrm{#1}}}

\newcommand{\vect}[1]{\vec{#1}}
\newcommand{\nhat}{\ensuremath{\hat{n}}}

\newcommand{\vecprime}[1]{\vect{#1}^{\,\prime}}

\renewcommand{\r}{\left(\vect{r}\right)}
\newcommand{\rp}{\left(\vecprime{r}\right)}
\newcommand{\rrp}{\left(\vect{r}, \vecprime{r}\right)}

\newcommand{\rrpk}[1][]{\left(k_{#1}, \vect{r}, \vecprime{r}\right)}
\newcommand{\rrpki}[1][]{\left(k_i, \vect{r}, \vecprime{r}\right)}

\newcommand{\matr}[1]{\bs{#1}}

\newcommand{\pvint}{\dashint}

\newcommand{\junk}[1] {}

\def\Xint#1{\mathchoice
{\XXint\displaystyle\textstyle{#1}}%
{\XXint\textstyle\scriptstyle{#1}}%
{\XXint\scriptstyle\scriptscriptstyle{#1}}%
{\XXint\scriptscriptstyle\scriptscriptstyle{#1}}%
\!\int}
\def\XXint#1#2#3{{\setbox0=\hbox{$#1{#2#3}{\int}$}
\vcenter{\hbox{$#2#3$}}\kern-.5\wd0}}

\def\dashint{\Xint-}

\newcommand*\widebar[1]{%
  \hbox{%
    \vbox{%
      \hrule height 0.5pt 
      \kern0.3ex
      \hbox{%
        \kern-0.05em
        \ensuremath{#1}%
        \kern-0.05em
      }%
    }%
  }%
} 

\newcommand{\mathsout}[1]
{\bgroup\mathchoice
  {\sbox0{$\displaystyle{#1}$}%
    \usebox0\hspace{-\wd0}%
    \rule[0.5\ht0-0.5\dp0-.5pt]{.7\wd0}{1pt}\hspace{.3\wd0}}%
  {\sbox0{$\textstyle{#1}$}%
    \usebox0\hspace{-\wd0}%
    \rule[0.5\ht0-0.5\dp0-.5pt]{.7\wd0}{1pt}\hspace{.3\wd0}}%
  {\sbox0{$\scriptstyle{#1}$}%
    \usebox0\hspace{-\wd0}%
    \rule[0.5\ht0-0.5\dp0-.5pt]{.7\wd0}{1pt}\hspace{.3\wd0}}%
  {\sbox0{$\scriptscriptstyle{#1}$}%
    \usebox0\hspace{-\wd0}%
    \rule[0.5\ht0-0.5\dp0-.5pt]{.7\wd0}{1pt}\hspace{.3\wd0}}%
\egroup}

\renewcommand{\epsilon}{\varepsilon}

\newcommand{\opL}{\ensuremath{\mathcal{L}}} 
\newcommand{\opK}{\ensuremath{\mathcal{K}}} 
\newcommand{\opKpv}{\mathsout{\opK}} 

\newcommand{\LmatA}[1][]{{\matr{L}^{(A)}_{#1}}}
\newcommand{\LmatPhi}[1][]{{\matr{L}^{(\phi)}_{#1}}}
\newcommand{\Kmat}[1][]{{\matr{K}_{#1}}}
\newcommand{\Kpvmat}[1][]{{\mathsout{\matr{K}}_{#1}}}

\newcommand{\Pxout}{\ensuremath{\matr{I}_{\times}}} 

\newcommand{\Pxin}{\ensuremath{\matr{I}_{\times}}} 

\newcommand{\figref}[1]{Fig.~\ref{#1}}
\newcommand{\secref}[1]{Section~\ref{#1}}
\newcommand{\tabref}[1]{Table~\ref{#1}}

\usepackage{tcolorbox}

\usepackage{textpos}

\newcommand{\mySubtitle}[1]%
{%
	\begin{textblock}{14.0}(0.7, 2.9)
		\textbf{#1}%
	\end{textblock}%
}%

\newcommand{\redcol}{black}

\newcommand{\red}[1]{\textcolor{\redcol}{#1}}
\newcommand{\green}[1]{\textcolor{green!0!black}{#1}}

\newcommand{\Ecolor}[1]{{#1}}
\newcommand{\Hcolor}[1]{{#1}}
\newcommand{\Jcolor}[1]{{#1}}

\newcommand{\Er}[1][]{\Ecolor{\vect{E}_{#1}\r}}

\newcommand{\Etr}[1][]{\Ecolor{\nhat_{#1} \times \vect{E}_{#1}\r}}
\newcommand{\Etrp}[1][]{\Ecolor{\nhat_{#1}' \times \vect{E}_{#1}\rp}}

\newcommand{\Emat}[1][]{\Ecolor{\matr{E}_{#1}}}

\newcommand{\Hr}[1][]{\Hcolor{\vect{H}_{#1}\r}}

\newcommand{\Htr}[1][]{\Hcolor{\nhat_{#1} \times \vect{H}_{#1}\r}}
\newcommand{\Htrp}[1][]{\Hcolor{\nhat_{#1}' \times \vect{H}_{#1}\rp}}
\newcommand{\Hmat}[1][]{\Hcolor{\matr{H}_{#1}}}

\newcommand{\Jr}[1][]{\ensuremath{\Jcolor{\vect{J}_{#1}\r}}}

\newcommand{\Jmat}[1][]{\ensuremath{\Jcolor{\matr{J}_{#1}}}}

\newcommand{\Grrp}[1][]{\ensuremath{\green{G_{#1}\rrp}}}
\newcommand{\Grrpk}[1][]{\ensuremath{\green{G_{#1}\rrpk}}}

\newcommand{\Grrpki}[1][]{\ensuremath{\green{G\rrpki}}}

\newcommand{\rhor}[1][]{\red{\rho_{#1}\r}}
\newcommand{\rhorp}[1][]{\red{\rho_{#1}\rp}}
\newcommand{\rhomat}[1][]{\red{\matr{\rho}_{#1}}}

\makeatletter
\newlength\numerator@height
\newlength\frac@height
\newsavebox\numerator@box
\newsavebox\frac@box

\newcommand\dfracparens[3]{%
	\sbox{\numerator@box}{\ensuremath{#1}}%
	\sbox{\frac@box}{\ensuremath{\dfrac{#1}{#2}}}%
	\settoheight{\frac@height}{\usebox{\frac@box}}%
	\settoheight{\numerator@height}{\usebox{\numerator@box}}%
	\addtolength{\frac@height}{-\numerator@height}%
	\usebox{\frac@box}%
	\raisebox{\frac@height}{%
		\( \left( {#3} \right)
		\)}%
}
\makeatother

%% file: introduction.tex
\section{Introduction}

\IEEEPARstart{T}{he} full-wave electromagnetic simulation of penetrable objects is crucial in a wide range of applications, spanning multiple scales of frequency, object size, and material properties.
For example, metamaterials and metasurfaces are often multiple wavelengths in size, but may contain sub-wavelength unit cells with intricate geometries.
Furthermore, the unit cells may be made of conductive~\cite{SRRmetaCond} or dielectric~\cite{optmeta01} materials embedded in a layered substrate.
Therefore, for an accurate simulation, the fields must be modeled both inside and outside each unit cell, while taking into account the layered surrounding medium.
Other applications include the design and analysis of ground-penetrating radar systems for the detection and reconstruction of objects buried beneath layers of soil and rock~\cite{gpr01_part1,gpr02,convcorr_AIM}.
The objects may be composed of dielectrics, conductors, or both, and may vary greatly in size and material properties.
Multiscale penetrable structures are also encountered in electronic packages and integrated circuits, where metallic traces and dielectric inclusions may occur in close proximity to each other.
Examples include air-backed or suspended on-chip inductor coils~\cite{indair01,indair02,indair03}.

The boundary element method (BEM), which is based on a surface integral representation of Maxwell's equations~\cite{book:colton}, has emerged as an important technique for full-wave electromagnetic modeling, since it allows three-dimensional problems to be modeled with a two-dimensional surface mesh~\cite{ChewIEM,gibson}.
For penetrable objects, the BEM requires solving an internal problem to model fields within objects, and an external problem to model the coupling between them.

Several BEM formulations have been proposed for modeling penetrable objects, but the existing methods have limitations, particularly for multiscale structures embedded in layered media.
The Poggio-Miller-Chang-Harrington-Wu-Tsai (PMCHWT)~\cite{PMCHWT02,PMCHWT03,PMCHWT04} and related formulations are well-established, but are only accurate when the contrast in electrical properties between adjacent media is not large.
This places restrictions on the types of materials which can be modeled.

The generalized impedance boundary condition (GIBC)~\cite{gibc} formulation and related methods~\cite{agibc,gibcHmatDanJiao} can handle a large contrast in material properties, but require the vector double-layer potential operator~\cite{book:colton} (often referred to as the $\opK$ operator) for the external problem.
When the multilayer Green's function (MGF)~\cite{MGF01} is used to model layered surrounding media, the double-layer operator requires computing the curl of the MGF in addition to the MGF itself, which significantly increases the computational cost~\cite{AWPLSLIM} and code complexity.
Furthermore, GIBC-based methods involve inner products between Rao-Wilton-Glisson (RWG) basis functions~\cite{RWG} and rotated RWG functions, which leads to a singular identity operator in both the internal and external problems.
This operator can deteriorate the convergence of iterative solvers for structures involving both conductive and dielectric objects, as will be shown in~\secref{sec:results}.
The enhanced augmented electric field integral equation (eAEFIE) formulation was proposed more recently for dielectrics~\cite{eaefie01} and lossy conductors~\cite{eaefie02}, and extends the original augmented electric field integral equation (AEFIE)~\cite{aefie1,aefie2}, which was applicable only to perfect conductors.
Unlike the GIBC, the eAEFIE uses the expensive Buffa-Christiansen (BC)~\cite{BCorig} functions, defined on a dual mesh, to obtain excellent convergence of iterative solvers for multiscale structures.
However, the eAEFIE also requires the double-layer potential operator in the external problem.
Moreover, this double-layer operator involves the BC functions, which compounds the associated computational cost.

The differential surface admittance (DSA) approach~\cite{DSA01,DSA06,UTK_AWPL2017} was proposed to avoid the double-layer potential operator in the external problem, but leads to a poorly conditioned system matrix even for relatively simple structures, as shown in~\cite{APS2020,AWPLSLIM}, in part because it also involves a singular identity operator.
Even when combined with BC functions~\cite{APS2019}, the DSA approach still involves two hypersingular operators which suffer from low-frequency breakdown~\cite{lfbreakdown}.
A formulation that combines the DSA approach with the Calder\'on preconditioner~\cite{BCcalderon} was proposed to obtain a well-conditioned matrix~\cite{DSA_Calderon_HDC}, which also uses BC functions defined on a dual mesh.
However, this formulation includes at least one matrix operator which involves both the MGF and BC functions, which significantly increases the computational cost associated with that matrix.
The same is true for the single-source formulation based on the Calder\'on preconditioner presented in~\cite{calderonSS}.
Well conditioned methods based on Helmholtz projectors, which avoid the use of BC functions, have also been proposed~\cite{RFCMP01,RFCMP02,RFCMP03}, but only for perfect electric conductors.

More recently, the single-layer impedance matrix (SLIM)~\cite{AWPLSLIM} formulation was proposed for lossy conductors.
The SLIM formulation, like the DSA, avoids the double-layer operator in the external problem, but is well conditioned without the need for BC functions.
This makes the SLIM approach well suited for modeling conductors embedded in layered media, such as high-speed electrical interconnects, and amenable to the use of iterative solvers~\cite{EPEPS2020}.
However, the SLIM formulation was proposed primarily for modeling lossy conductors, and may lead to a poorly conditioned system matrix for multiscale structures involving both dielectric and conductive objects.
Also, some operators involved in the SLIM approach also suffer from low-frequency breakdown.

In this article, we propose a new formulation for modeling penetrable objects in layered media, which combines the benefits of the eAEFIE and SLIM formulations, while avoiding their drawbacks.
The formulation we propose is suitable for modeling both dielectric and conductive objects, and accurately captures the frequency-dependent variations of skin depth in the latter.
The proposed technique avoids the double-layer operator in the external problem, making it more computationally efficient than the eAEFIE for objects embedded in layered media.
Unlike the SLIM approach, the proposed single-layer formulation also eliminates the low-frequency breakdown by using the AEFIE in both the internal and external problems.
The proposed method employs BC functions defined on a dual mesh to obtain excellent conditioning of the system matrix, even for dielectrics and multiscale structures, unlike the SLIM approach.
An important and distinctive feature of the proposed formulation is that integral operators involving the expensive BC functions only appear in the internal problem, which has two advantages:
\begin{itemize}
	\item All operators which involve BC functions require only the simple homogeneous Green's function, which is available in closed form, rather than the MGF, which is expensive to compute and not expressible in closed form.
	\item All interactions involving the BC functions are local to each object, which is particularly advantageous when the same object occurs in multiple locations.
\end{itemize} 
Therefore, for layered background media, the proposed formulation leverages the benefits of BC functions but greatly alleviates the associated computational cost, compared to the eAEFIE and the method in~\cite{DSA_Calderon_HDC}.

The paper is organized as follows: \secref{sec:formulation} describes the proposed formulation and provides a detailed comparison to existing techniques. \secref{sec:acc} summarizes how acceleration algorithms are incorporated into the proposed method to handle large problems efficiently. In \secref{sec:results}, we demonstrate the efficiency and good conditioning of the proposed method, compared to existing formulations, for multiscale numerical examples drawn from different applications.
Concluding remarks are provided in \secref{sec:conclusion}.


%% file: methods.tex
\section{Proposed Formulation}\label{sec:formulation}

\begin{figure*}[t]
	\centering
	\subfloat[][]
	{\includegraphics[width=0.4\linewidth]{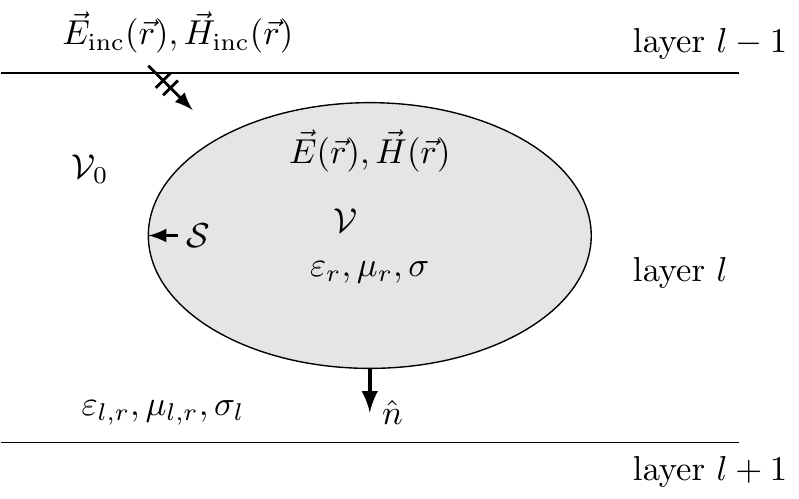}\label{fig:geom:a}}\hspace{6mm}
	\subfloat[][]
	{\includegraphics[width=0.4\linewidth]{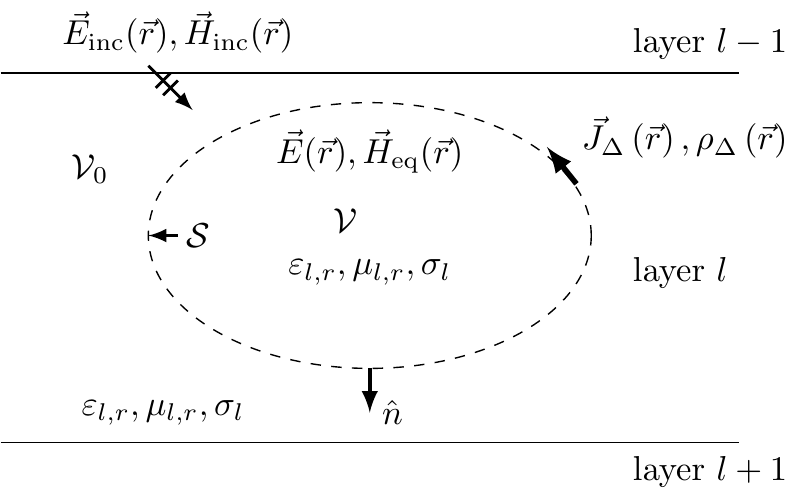}\label{fig:geom:b}}
	\caption{Geometric setup considered in this article, showing (a) the original configuration, and (b) the equivalent configuration.}\label{fig:geom}
\end{figure*}

We consider the problem of electromagnetic scattering from a penetrable object occupying volume $\mathcal{V}$, bounded by a surface $\mathcal{S}$ with outward unit normal vector $\nhat$, as shown in \figref{fig:geom:a}.
The object is assumed to be homogeneous with relative permittivity $\epsilon_r$, relative permeability $\mu_r$, and conductivity $\sigma$.
The object may be embedded in free space, or in layer $l$ of a stratified background medium.
The general case of multiple objects is discussed in \secref{sec:multiobject}.
The background medium is denoted by $\mathcal{V}_0$.
In the case of a layered medium, the $l^{\mathrm{th}}$ layer has relative permittivity $\epsilon_{l,r}$, relative permeability $\mu_{l,r}$, and conductivity $\sigma_l$.
The time-harmonic fields ${[\Er[\mathrm{inc}], \Hr[\mathrm{inc}]]}$, ${\vect{r} \in \mathcal{V}_0}$, with cyclical frequency $\omega$, are incident on $\mathcal{S}$.
This leads to the field distributions ${[\Er, \Hr]}$, ${\vect{r} \in \mathcal{V}}$.

\subsection{Overview}\label{sec:overview}
We first provide a qualitative overview of the proposed formulation, followed by the full mathematical derivation.
In most existing BEM formulations for penetrable objects, such as the GIBC~\cite{gibc} and the eAEFIE~\cite{eaefie01}, surface integral equations~(SIEs) are formulated for the regions external and internal to the object, and solved jointly to compute the fields tangential to~$\mathcal{S}$.
This necessitates the use of the vector double-layer potential operator associated with~$\mathcal{V}_0$.
Instead, in the proposed method, we introduce additional SIEs which are associated with an equivalent configuration, where the object is replaced by the surrounding medium and an equivalent electric surface current density~\cite{DSA01}, as shown in \figref{fig:geom:b}.
The additional degrees of freedom introduced via the equivalent configuration allow deriving a final system of equations where the double-layer operator associated with~$\mathcal{V}_0$ does not appear.
Unlike other methods which use the equivalence principle in a similar manner~\cite{UTK_AWPL2017,DSA06,DSA08,AWPLSLIM}, we use BC basis functions~\cite{BCorig} to express the electric field tangential to~$\mathcal{S}$, so that the basis and testing functions for each integral operator reside in the appropriate function space~\cite{APmagUltimate,eaefie01}.
This leads to a system of equations which is well conditioned over multiple scales of frequency and material parameters.
Furthermore, we introduce the electric charge density as an unknown~\cite{aefie1} in both the original and equivalent scenarios, unlike existing methods, to ensure stability at low frequencies.
Together, these strategies allow us to devise a final system of equations which is as well conditioned as state-of-the-art formulations~\cite{eaefie01}, but is significantly more efficient when~$\mathcal{V}_0$ consists of a stratified substrate.

\subsection{Internal Problem}\label{sec:internal}

First, we derive a set of surface integral equations to model field distributions interior to~$\mathcal{V}$, which we refer to as the internal problem.
The tangential electric and magnetic fields on $\mathcal{S}$ can be related via the augmented electric field integral equation~(AEFIE) formulation~\cite{aefie1,aefie2}.
The original AEFIE in~\cite{aefie1} was developed for perfect conductors, and must be modified for penetrable objects to include the tangential electric field on~$\mathcal{S}$~\cite{eaefie01}, as
\begin{subequations}
\begin{multline}
	\mu_r\,\nhat \times \opL^{(A)} \bigl[jk_0\,\Htrp\bigr]
	- \epsilon_{c,r}^{-1}\,\,\nhat \times \nabla \opL^{(\phi)} \bigl[c_0\,\rhorp\bigr]\\
	- \nhat \times \opKpv \bigl[\eta_0^{-1}\,\Etrp\bigr]
	- \frac{1}{2} \Bigl(\eta_0^{-1}\,\Etr\Bigr) = 0,\label{eq:AEFIEint1}
\end{multline}
\vspace{-5mm}
\begin{align}
	\nabla\cdot\Bigl(jk_0\,\Htr\Bigr) - k_0^2\Bigl(c_0\,\rhor\Bigr) = 0,\label{eq:AEFIEint2}
\end{align}\label{eq:AEFIEint}%
\end{subequations}%
where $\vect{r}, \vect{r}\,' \in \mathcal{S}$, and $c_0$, $k_0$ and $\eta_0$ are the speed of light, wave number and wave impedance in free space.
Primed coordinates represent source points, while unprimed coordinates represent observation points.
The relative complex permittivity $\epsilon_{c,r}$ is defined as
\begin{align}
	\epsilon_{c,r} = \epsilon_r - j\frac{\sigma}{\omega\epsilon_0}.\label{eq:epsrc}
\end{align}
Equation~\eqref{eq:AEFIEint1} is the conventional electric field integral equation~(EFIE)~\cite{ChewWAF}, but with charge density~$\rhor$ taken as a separate unknown~\cite{aefie2}, while~\eqref{eq:AEFIEint2} is the continuity equation relating~$\rhor$ to~$\Htr$.
The integral operators in~\eqref{eq:AEFIEint1} are defined as~\cite{ChewIEM}
\begin{align}
	\opL^{(A)}\bigl[\vect{X}\rp\bigr] &= \int_\mathcal{S} d\mathcal{S}'\, \vect{X}\rp\,\Grrpk,\label{eq:opLAhgf}\\
	\opL^{(\phi)}\bigl[a\rp\bigr] &= \int_\mathcal{S} d\mathcal{S}'\, a\rp\,\Grrpk,\label{eq:opLphihgf}\\
	\opKpv\bigl[\vect{X}\rp\bigr] &= \pvint_\mathcal{S} d\mathcal{S}'\, \nabla\Grrpk \times \vect{X}\rp,\label{eq:opKhgf}
\end{align}
where $k = \sqrt{\epsilon_{c,r}\mu_r}\,k_0$ is the wave number associated with the object's material.
The homogeneous Green's function is
\begin{align}
  \Grrp = \dfrac{e^{-j k r}}{4\pi r},\label{eq:hgf}
\end{align}
where $r = |\vect{r} - \vect{r}^{\,'}|$.

A triangular mesh is generated for the surface of the object.
In~\eqref{eq:AEFIEint}, $\Htr$ is expanded with divergence-conforming area-normalized $\text{RWG}$ functions~\cite{RWG}, while $\Etr$ is expanded with Buffa Christiansen~(BC) functions~\cite{BCorig} defined on a barycentric refinement of the original mesh.
Quantity~$\rhor$ is expanded with area-normalized pulse functions.
Equation~\eqref{eq:AEFIEint1} is tested with ${\nhat \times \text{RWG}}$ functions, while~\eqref{eq:AEFIEint2} is tested with pulse functions.
The gradient operator in~\eqref{eq:AEFIEint1} is transferred to the testing function~\cite{gibson}, and the discretized AEFIE is
\begin{multline}
	{%
	\begin{bmatrix}
		\mu_r\, \LmatA & -\epsilon_{c,r}^{-1}\,\matr{D}^T\LmatPhi\matr{B} \\
		\matr{F}\matr{D} & -k_0^2\,\matr{I}
	\end{bmatrix}}
	\begin{bmatrix}
		jk_0\Hmat \\ c_0\rhomat
	\end{bmatrix} \\=
	{%
	\begin{bmatrix}
		\left( \Kpvmat + \frac{1}{2}\,\Pxin \right) & \matr{0} \\
		\matr{0} & \matr{I}
	\end{bmatrix}}
	\begin{bmatrix}
		\eta_0^{-1}\,\Emat \\ \matr{0}
	\end{bmatrix},\label{eq:AEFIEintdis}
\end{multline}
where $\LmatA$, $\LmatPhi$ and $\Kpvmat$ are the discretized $\opL^{(A)}$, $\opL^{(\phi)}$ and $\opKpv$ operators, respectively.
Matrix~$\matr{D}$ is an incidence matrix, while~$\matr{F}$ and~$\matr{B}$ are used to enforce charge neutrality~\cite{aefie2}.
Definitions of~$\matr{D}$,~$\matr{F}$ and~$\matr{B}$ can be found in~\cite{aefie2}.
Matrix $\matr{I}$ is the identity, while $\Pxout$ is the well-conditioned Gram matrix obtained by testing $\text{BC}$ basis functions with ${\nhat\times\text{RWG}}$ functions~\cite{BCcalderon}.
Column vectors $\Hmat$, $\rhomat$ and $\Emat$ contain the coefficients of the basis functions associated with $\Htr$, $\rhor$, and $\Etr$, respectively.
For highly conductive objects, specialized integration routines~\cite{gibc} are used to compute the entries of $\LmatA$, $\LmatPhi$ and $\matr{K}$, to capture the fast oscillations of the Green's function.
%

Next, in order to avoid the double-layer operator in the external problem, the surface equivalence principle~\cite{EMharrington} is invoked to replace the object with the surrounding material, while keeping $\Etr$ unchanged for ${\vec{r}\in\mathcal{S}}$~\cite{DSA01}.
A pair of equivalent differential sources, $\Jr[\Delta]$ and $\rhor[\Delta]$,
\begin{align}
	\begin{bmatrix}
  		\Jr[\Delta] \\ \rhor[\Delta]
  	\end{bmatrix} = 
  	\begin{bmatrix}
  		\Htr - \Hcolor{\nhat\, \times\, }\Hr[\mathrm{eq}]\\
  		\rhor - \rhor[\mathrm{eq}]
  	\end{bmatrix}
  	 ,\label{eq:diffsources}
\end{align}
must be introduced on $\mathcal{S}$ to keep fields in $\mathcal{V}_0$ unchanged, as shown in \figref{fig:geom:b}.
In \eqref{eq:diffsources}, ${\Hcolor{\nhat\,\times}\Hr[\mathrm{eq}]}$ and $\rhor[\mathrm{eq}]$ are the tangential magnetic field and surface charge density on $\mathcal{S}$, respectively, in the equivalent configuration.
Keeping $\Etr$ unchanged is the key to avoiding the double-layer potential operator in the external problem~\cite{DSA01}.

In~\eqref{eq:diffsources}, $\vect{J}_{\Delta}$, $\Htr$ and ${\Hcolor{\nhat\, \times\, }\Hr[\mathrm{eq}]}$ are expanded with $\text{RWG}$ functions, while $\rhor[\Delta]$, $\rhor$ and $\rhor[\mathrm{eq}]$ are expanded with pulse functions. This yields
\begin{align}
	\begin{bmatrix}
  		\Jmat[\Delta] \\ \rhomat[\Delta]
  	\end{bmatrix} =
  	\begin{bmatrix}
  		\Hmat - \Hmat[\mathrm{eq}] \\ \rhomat - \rhomat[\mathrm{eq}]
  	\end{bmatrix},\label{eq:diffsourcesdis}
\end{align}
where $\Jmat[\Delta]$ and $\Hmat[\mathrm{eq}]$ contain the coefficients of the basis functions associated with $\vect{J}_{\Delta}$ and ${\Hcolor{\nhat\, \times\, }\Hr[\mathrm{eq}]}$, respectively.
Column vectors $\rhomat[\Delta]$ and $\rhomat[\mathrm{eq}]$ contain the basis function coefficients associated with $\rhor[\Delta]$ and $\rhor[\mathrm{eq}]$, respectively.

The AEFIE~\cite{aefie2} is invoked again in the equivalent configuration to relate ${\Hcolor{\nhat\, \times\, }\Hr[\mathrm{eq}]}$, $\Etr$ and $\rhor[\mathrm{eq}]$, which in discrete form is
\begin{multline}
	{%
		\begin{bmatrix}
			\mu_{l,r}\, \LmatA[l] & -\epsilon_{l,r}^{-1}\,\matr{D}^T\LmatPhi[l]\matr{B} \\
			\matr{F}\matr{D} & -k_0^2\,\matr{I}
	\end{bmatrix}}
	\begin{bmatrix}
		jk_0\Hmat[\mathrm{eq}] \\ c_0\rhomat[\mathrm{eq}]
	\end{bmatrix} \\=
	{%
		\begin{bmatrix}
			\left( \Kpvmat[l] + \frac{1}{2}\,\Pxin \right) & \matr{0} \\
			\matr{0} & \matr{I}
	\end{bmatrix}}
	\begin{bmatrix}
		\eta_0^{-1}\,\Emat \\ \matr{0}
	\end{bmatrix}.\label{eq:AEFIEeqvdis}
\end{multline}
Equation~\eqref{eq:AEFIEeqvdis} is analogous to~\eqref{eq:AEFIEintdis}, except that it corresponds to the equivalent configuration where the objects have been replaced by the surrounding medium~\cite{DSA01}.
To obtain~\eqref{eq:AEFIEeqvdis},~$\Etr$ was expanded with $\text{BC}$ functions, and the AEFIE was tested with~${\nhat\times\text{RWG}}$ and pulse functions, as in the original configuration.
In~\eqref{eq:AEFIEeqvdis},~$\LmatA[l]$,~$\LmatPhi[l]$ and~$\Kpvmat[l]$ are the discretized~$\opL^{(A)}$,~$\opL^{(\phi)}$, and~$\opKpv$ operators, respectively, involving the homogeneous Green's function associated to layer $l$ of the background medium.
In order to eliminate~$\Hmat[\mathrm{eq}]$ and~$\rhomat[\mathrm{eq}]$ from the final system of equations,~\eqref{eq:AEFIEeqvdis} is rearranged as
\begin{align}
	\begin{bmatrix}
		jk_0\Hmat[\mathrm{eq}] \\ c_0\rhomat[\mathrm{eq}]
	\end{bmatrix} =
	{%
		\begin{bmatrix}
			\matr{Y}_{11} & \matr{Y}_{12} \\
			\matr{Y}_{21} & \matr{Y}_{22}
	\end{bmatrix}}
	\begin{bmatrix}
		\eta_0^{-1}\,\Emat \\ \matr{0}
	\end{bmatrix} =
	\begin{bmatrix}
		\eta_0^{-1}\,\matr{Y}_{11}\Emat \\ \eta_0^{-1}\,\matr{Y}_{21}\Emat
	\end{bmatrix},\label{eq:invAEFIEeqvdis}
\end{align}
where we have introduced the admittance-like matrix operator,
\begin{multline}
	{%
	\begin{bmatrix}
		\matr{Y}_{11} & \matr{Y}_{12} \\
		\matr{Y}_{21} & \matr{Y}_{22}
	\end{bmatrix}} =
	{%
	\setlength\arraycolsep{2pt}
	\begin{bmatrix}
		\mu_{l,r}\, \LmatA[l] & -\epsilon_{l,r}^{-1}\,\matr{D}^T\LmatPhi[l]\matr{B} \\
		\matr{F}\matr{D} & -k_0^2\,\matr{I}
	\end{bmatrix}^{-1}}\cdot\\
	{%
	\setlength\arraycolsep{2pt}
	\begin{bmatrix}
		\left( \Kpvmat[l] + \frac{1}{2}\,\Pxin \right) & \matr{0} \\
		\matr{0} & \matr{I}
	\end{bmatrix}}.\label{eq:Yeqdef}
\end{multline}
In~\secref{sec:acc}, we will describe how the matrix inversion in~\eqref{eq:Yeqdef} can be handled efficiently without the need for direct factorization, even for large problems. 
\secref{sec:results} demonstrates that the proposed method compares favorably to existing techniques such as the eAEFIE~\cite{eaefie01} for layered medium problems, despite the additional matrix operators required in the equivalent configuration~\eqref{eq:AEFIEeqvdis}.

\subsection{External Problem}\label{sec:external}

Having applied the equivalence theorem, the objects have been replaced by the differential sources~$\Jr[\Delta]$ and~$\rhor[\Delta]$, which radiate into the surrounding homogeneous or layered medium.
These sources are used to model the external problem by relating them to the incident fields via the AEFIE~\cite{aefie2}.
Expanding~$\Jr[\Delta]$,~$\Etr$ and~$\rhor[\Delta]$ with RWG, BC and pulse functions, respectively, then testing the AEFIE as before,
\begin{multline}
	{%
	\begin{bmatrix}
		\LmatA[m] & -\matr{D}^T\LmatPhi[m]\matr{B} \\
		\matr{F}\matr{D} & -k_0^2\,\matr{I}
	\end{bmatrix}}
	\begin{bmatrix}
		jk_0\Jmat[\Delta] \\ c_0\rhomat[\Delta]
	\end{bmatrix} \\+
	{%
	\begin{bmatrix}
		\Pxin & \matr{0} \\
		\matr{0} & \matr{I}
	\end{bmatrix}}
	\begin{bmatrix}
		\eta_0^{-1}\,\Emat \\ \matr{0}
	\end{bmatrix} =
	\begin{bmatrix}
		\eta_0^{-1}\,\Emat[\mathrm{inc}] \\ \matr{0}
	\end{bmatrix},\label{eq:AEFIEextdis}
\end{multline}
where subscript~$\mathrm{m}$ indicates the use of the dyadic multilayer Green's function~(MGF) of the background medium as the kernel of the associated integral operator~\cite{MGF01}.

\subsection{Final System Matrix}\label{sec:finalsys}

We now seek a final system of equations in terms of the unknowns~$\Hmat$,~$\Emat$ and~$\rhomat$.
Using~\eqref{eq:diffsourcesdis} in~\eqref{eq:AEFIEextdis} and rearranging, we obtain
\begin{multline}
	{%
	\begin{bmatrix}
		\LmatA[m] & -\matr{D}^T\LmatPhi[m]\matr{B} \\
		\matr{F}\matr{D} & -k_0^2\,\matr{I}
	\end{bmatrix}}
	\begin{bmatrix}
  		jk_0\,\Hmat \\ c_0\,\rhomat
  	\end{bmatrix} \\-
  	{%
	\begin{bmatrix}
		\LmatA[m] & -\matr{D}^T\LmatPhi[m]\matr{B} \\
		\matr{F}\matr{D} & -k_0^2\,\matr{I}
	\end{bmatrix}}
	\begin{bmatrix}
  		jk_0\,\Hmat[\mathrm{eq}] \\ c_0\,\rhomat[\mathrm{eq}]
  	\end{bmatrix}
  	\\+
	{%
	\begin{bmatrix}
		\Pxin & \matr{0} \\
		\matr{0} & \matr{I}
	\end{bmatrix}}
	\begin{bmatrix}
		\eta_0^{-1}\,\Emat \\ \matr{0}
	\end{bmatrix} =
	\begin{bmatrix}
		\eta_0^{-1}\,\Emat[\mathrm{inc}] \\ \matr{0}
	\end{bmatrix}.\label{eq:AEFIEextdisY}
\end{multline}
Next, using~\eqref{eq:invAEFIEeqvdis} in~\eqref{eq:AEFIEextdis} and rearranging terms gives
\begin{align}
	{%
	\begin{bmatrix}
		\LmatA[m] & \matr{C}_1 & -\matr{D}^T\LmatPhi[m]\matr{B} \\
		\matr{F}\matr{D} & \matr{C}_2 & -k_0^2\,\matr{I}
	\end{bmatrix}}
	\begin{bmatrix}
  		jk_0\,\Hmat \\ \eta_0^{-1}\,\Emat \\ c_0\,\rhomat
  	\end{bmatrix}
	=
	\begin{bmatrix}
		\eta_0^{-1}\,\Emat[\mathrm{inc}] \\ \matr{0}
	\end{bmatrix},\label{eq:prefinalsystem}
\end{align}
where
\begin{align}
	\matr{C}_1 &= \Pxin - \left(\LmatA[m] \matr{Y}_{11} - \matr{D}^T\LmatPhi[l]\matr{B} \matr{Y}_{21}\right),\label{eq:C1}\\
	\matr{C}_2 &= -\left(\matr{F}\matr{D} \matr{Y}_{11} -k_0^2\,\matr{Y}_{21}\right).\label{eq:C2}
\end{align}
Finally, we use the first equation of~\eqref{eq:AEFIEintdis} in~\eqref{eq:prefinalsystem} to obtain the proposed single-layer dual-mesh formulation,
\begin{align}
	{%
	\setlength\arraycolsep{4.5pt}
	\begin{bmatrix}
		\LmatA[m] & \matr{C}_1 & \matr{D}^T\LmatPhi[m]\matr{B} \\
		\mu_r\, \LmatA & \Kmat & \epsilon_{c,r}^{-1}\,\matr{D}^T\LmatPhi\matr{B} \\
		\matr{F}\matr{D} & \matr{C}_2 & k_0^2\,\matr{I}
	\end{bmatrix}}
	\begin{bmatrix}
  		jk_0\,\Hmat \\ \eta_0^{-1}\,\Emat \\ -c_0\,\rhomat
  	\end{bmatrix}
	=
	\begin{bmatrix}
		\frac{\Emat[\mathrm{inc}]}{\eta_0} \\ \matr{0} \\ \matr{0}
	\end{bmatrix},\label{eq:finalsystem}
\end{align}
where
\begin{align}
	\Kmat = \Kpvmat + \frac{1}{2}\Pxin.\label{eq:Kmatdef}
\end{align}

Although the double-layer potential operator for the external region is avoided, double-layer operators associated with the internal problem still appear in~\eqref{eq:finalsystem}, and may be poorly conditioned at very low frequencies for multiply-connected geometries~\cite{KNULLSPACE01}\footnote{We are grateful to an anonymous reviewer for pointing this out.}.
This occurs due to an approximate null space of the double-layer operator in the static limit, with a dimension equal to the genus of the associated object~\cite{KNULLSPACE01}.
Accurate results at very low frequencies can still be obtained by using an alternative discretization and testing scheme~\cite{KNULLSPACE02}, by using Calder\'on-like formulations involving quasi-Helmholtz projectors~\cite{KNULLSPACE03}, or by imposing an additional consistency condition~\cite{KNULLCONSIST}.
In this work, we do not address the low-frequency null space of the double-layer operator.
However, the numerical example in \secref{sec:results:ind} shows that accurate results can still be obtained for sub-wavelength structures of practical relevance with a genus larger than $0$.

The system matrices in both~\eqref{eq:AEFIEeqvdis} and~\eqref{eq:finalsystem} represent a generalized saddle point problem~\cite{saddle,aefie1,eaefie01}.
In the static limit, the system matrix in~\eqref{eq:finalsystem} scales with frequency as
\begin{align}
	{%
	\setlength\arraycolsep{4.5pt}
	\begin{bmatrix}
		\mathcal{O}\left(1\right) & \mathcal{O}\left(1\right) & \mathcal{O}\left(1\right) \\
		\mu_r\, \mathcal{O}\left(1\right) & \mathcal{O}\left(1\right) & \epsilon_{c,r}^{-1}\,\mathcal{O}\left(1\right) \\
		\mathcal{O}\left(1\right) & \mathcal{O}\left(1\right) & 0
	\end{bmatrix}}\label{eq:DCscalediel}
\end{align}
for dielectrics, and
\begin{align}
	{%
	\setlength\arraycolsep{4.5pt}
	\begin{bmatrix}
		\mathcal{O}\left(1\right) & \mathcal{O}\left(1\right) & \mathcal{O}\left(1\right) \\
		\mu_r\, \mathcal{O}\left(1\right) & \mathcal{O}\left(1\right) & 0 \\
		\mathcal{O}\left(1\right) & \mathcal{O}\left(1\right) & 0
	\end{bmatrix}}\label{eq:DCscalecond}
\end{align}
for conductors.
The enforcement of charge neutrality via~$\matr{B}$ and~$\matr{F}$ ensures that the system matrix in~\eqref{eq:finalsystem} has full rank and is amenable to the use of an iterative solver at low frequencies, with the use of the preconditioner discussed in \secref{sec:acc}~\cite{aefie2,eaefie01}.
However, similar to the AEFIE~\cite{aefie2} and the eAEFIE~\cite{eaefie01}, the proposed method may yield an inaccurate solution for~$\Hmat$ for scattering and capacitive problems at very low frequencies, while the charge unknowns in~$\rhomat$ remain accurate because they dominate over~$jk_0\,\Hmat$~\cite{ChewAEFIEperturbOrig}.
This low frequency inaccuracy can be analyzed with a procedure similar to~\cite{ChewAEFIEperturbOrig}, where a perturbation approach is devised to remedy this issue.

\begin{figure}[t]
  \centering
  \includegraphics[width=\linewidth]{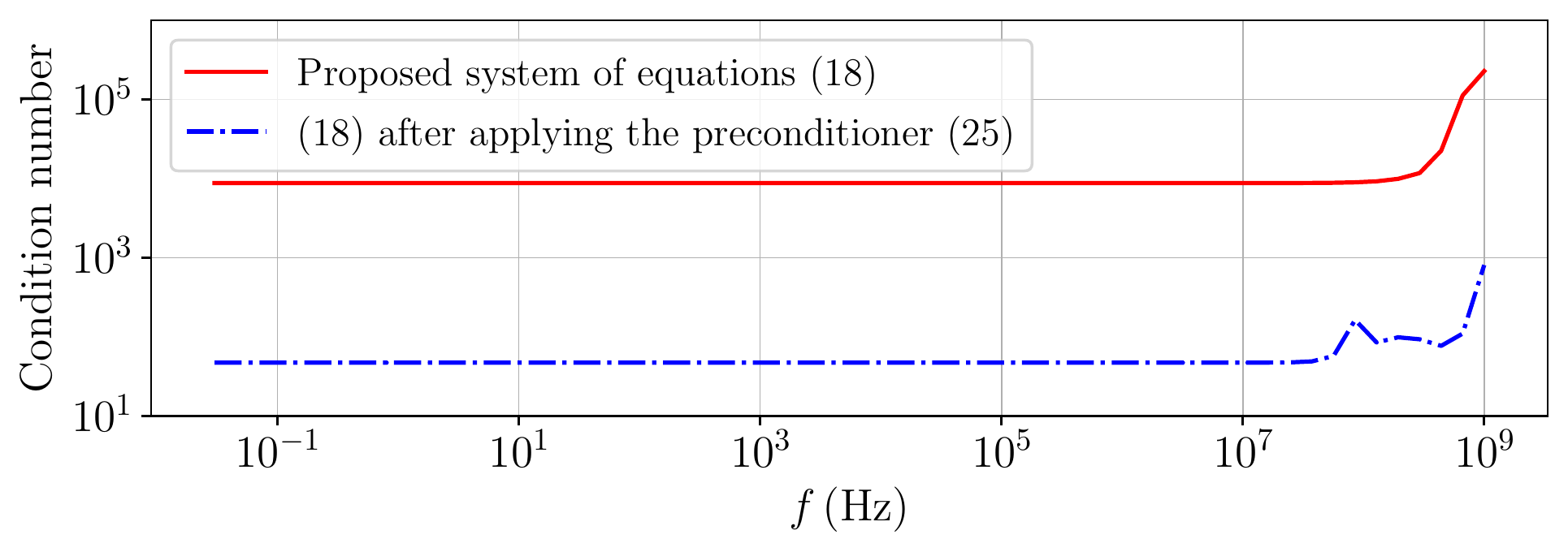}\\
  \caption{Matrix conditioning for the sphere in \secref{sec:finalsys}.}\label{fig:sphcond}
\end{figure}

To demonstrate numerically the stability of the proposed formulation at low frequencies, we consider a dielectric sphere in free space with diameter~$1\,$m and relative permittivity~$2.56$.
The sphere is meshed with~$752$ triangles, and the condition number of the system matrix in~\eqref{eq:finalsystem} is reported over a wide range of frequencies in \figref{fig:sphcond}.
Also reported is the condition number of the system matrix in~\eqref{eq:finalsystem} after applying the preconditioner~\eqref{eq:pcext} described in \secref{sec:acc}.
\figref{fig:sphcond} demonstrates that the system matrix associated with the proposed method remains stable at very low frequencies, and that the condition number can be significantly improved with an appropriate preconditioner, which is discussed further in \secref{sec:acc}.

\subsection{Structures with Multiple Objects}\label{sec:multiobject}

Next, we generalize the proposed formulation~\eqref{eq:finalsystem} to the case of structures containing multiple objects, which is easily accomplished with a block-diagonal concatenation of the internal problem matrices.
We assume that the objects are disjoint; the case of objects in contact with each other (and the resulting formation of junctions) is not considered here, but can be incorporated easily with the use of an appropriate set of boundary conditions.

Suppose that~$\matr{A}$ represents any of the blocks in the matrices in~\eqref{eq:AEFIEintdis} or~\eqref{eq:AEFIEeqvdis}.
For multiple objects, we may write
\begin{align}
	\matr{A} = \mathrm{diag}
	\begin{bmatrix}
		\matr{A}_1 & \matr{A}_2 & \cdots & \matr{A}_{N_{\mathrm{obj}}}
	\end{bmatrix},\label{eq:Aobj}
\end{align}
where~$N_{\mathrm{obj}}$ is the total number of disjoint objects.
The appropriate material properties associated with the~$i^{\mathrm{th}}$ object must be used in $\matr{A}_i$.
A key feature of the proposed formulation is that the expensive BC functions are involved only in the matrix operators of~\eqref{eq:AEFIEintdis} and~\eqref{eq:AEFIEeqvdis}, and not in the external problem~\eqref{eq:AEFIEextdis}.
The only exception is the sparse Gram matrix~$\Pxin$ in~\eqref{eq:AEFIEextdis}, which does not involve a Green's function and is very quick to compute.
Therefore, thanks to the block-diagonal structure implied by~\eqref{eq:Aobj}, integral operators involving BC functions are local to each object, while the inter-object couplings encoded in~$\LmatA[\mathrm{m}]$ and~$\LmatPhi[\mathrm{m}]$ only involve RWG functions.
This is computationally advantageous for multi-object structures, because it implies fewer source-test interactions involving BC functions, compared to the eAEFIE~\cite{eaefie01} and the method in~\cite{DSA_Calderon_HDC}, where BC functions are involved in both intra- and inter-object interactions.
As shown in \secref{sec:results}, CPU time savings up to nearly~$7\times$ are achieved compared to the eAEFIE.

\subsection{Discussion and Comparison to the eAEFIE}\label{sec:comparison}

The form of~\eqref{eq:finalsystem} resembles that of the eAEFIE~\cite{eaefie01}, except for the blocks which contain~$\matr{C}_1$ and~$\matr{C}_2$.
The counterparts to~$\matr{C}_1$ and~$\matr{C}_2$ in the eAEFIE are
\begin{align}
	\matr{C}_1^{\text{\cite{eaefie01}}} &= \Kpvmat[\mathrm{m}] - \frac{1}{2}\Pxin,\\
	\matr{C}_2^{\text{\cite{eaefie01}}} &= \matr{0},
\end{align}
where~$\Kpvmat[\mathrm{m}]$ is associated to the external problem, and is the discretized double-layer operator~$\opKpv$ for the surrounding medium.
The key distinction of the proposed formulation is that we use the equivalence principle and the differential sources~$\Jr[\Delta]$ and~$\rhor[\Delta]$ to avoid~$\Kpvmat[\mathrm{m}]$.
This leads to a significant computational advantage over the eAEFIE in applications involving objects embedded in layered media, for the following reasons:
\begin{enumerate}
	\item\label{adv1} The proposed method requires only the single-layer operator in the external problem, which involves computing the MGF~\cite{MGF01}.
	In contrast, due to the presence of~$\Kpvmat[\mathrm{m}]$, the eAEFIE requires not only the MGF, but also its curl~\cite{MGF01}.
	The curl of the MGF is approximately as expensive to compute as the MGF itself, which can lead to a~$2\times$ increase in computational cost compared to a formulation that requires only the MGF~\cite{AWPLSLIM}. 
	\item\label{adv2} In the proposed formulation, BC functions only appear in integral operators associated with the \emph{internal} problem,~\eqref{eq:AEFIEintdis} and~\eqref{eq:AEFIEeqvdis}, which involve the simple homogeneous Green's function~\eqref{eq:hgf}.
	Therefore, the MGF and the BC functions never occur simultaneously in the same operator. 
	Instead in the eAEFIE,~$\Kpvmat[\mathrm{m}]$ is associated with the external problem, and involves not only the use of BC basis functions, but also the curl of the MGF.
	Considering point~\ref{adv1}) above, the presence of both the curl of the MGF \emph{and} BC functions in the \emph{same} operator is especially expensive, since computations involving BC functions can be~$6\times$ more expensive than those involving only RWG functions.
	This point contributes significantly to the computational advantage of the proposed approach compared to the eAEFIE, as verified in~\secref{sec:results}.
	\item\label{adv3} In the proposed formulation, the computational cost associated with BC functions can be further reduced for structures where the same object occurs in multiple locations.
	Examples include the use of repeated unit cells in antenna arrays and metasurfaces, and arrays of vias and ground bars often encountered in on-chip passive components and interconnects.
	In these cases, as long as the mesh is identical for each repeated object, the same internal problem matrices can be reused, leading to significant computational savings.
	Though this can also be done for the eAEFIE formulation, the reuse of internal problem matrices provides a greater advantage to the proposed formulation.
	This is because in the eAEFIE, BC functions are used in both the internal and external problems.
	Therefore, even in the best-case scenario of only one unique object in the structure, only some of the computations related to BC functions (those associated with the internal problem) can be reused.
	In the proposed method, BC functions only appear in the internal problem.
	This allows reusing \emph{all} computations associated with the BC functions, for the objects which are repeated.
	For simplicity, the internal problem matrices were not reused in any of the examples considered in \secref{sec:results}.
\end{enumerate}
%
At the same time, the proposed formulation~\eqref{eq:finalsystem} also retains three key features of the eAEFIE, which contribute to its good conditioning:
\begin{itemize}
	\item Each block along the diagonal is well conditioned, particularly~$\Kmat$, due the use of a dual mesh and BC functions.
	\item The blocks of~\eqref{eq:finalsystem} are strategically scaled similarly to the eAEFIE, which contributes towards stability at low frequencies.
	\item When using an iterative solver, the form of~\eqref{eq:finalsystem} allows using a powerful preconditioner similar to the one proposed for the eAEFIE, as described in \secref{sec:acc}, which significantly reduces the number of iterations required.
\end{itemize}
Therefore, the proposed method is able to exploit the features of the eAEFIE which lead to good conditioning, while offering a significant computational advantage for the case of layered surrounding media.

\section{Acceleration for Large Multiscale Problems}\label{sec:acc}

To handle large structures, we propose an acceleration scheme for the proposed formulation which provides an overall~${\sim}\mathcal{O}(N^{1.5}\log N)$ time complexity, where~$N$ is the total number of unknowns in~\eqref{eq:finalsystem}.
The adaptive integral method (AIM)~\cite{AIMbles} is employed in both the external and internal problems, along with an iterative solver.
In~\eqref{eq:finalsystem}, the matrices~$\LmatA[m]$ and~$\LmatPhi[m]$ are compressed via a specialization of the AIM for layered media~\cite{CPMT2019arxiv}.
Instead,~$\LmatA$,~$\LmatPhi$, and~$\Kmat$ are straightforward to compress via the conventional AIM for both dielectrics and conductors~\cite{TAPAIMin}, since they involve the homogeneous Green's function~\eqref{eq:hgf}.
In~\eqref{eq:C1} and~\eqref{eq:C2},~$\matr{Y}_{11}$ and~$\matr{Y}_{21}$ are not computed explicitly.
Instead, the products~$\matr{Y}_{11}\Emat^{(i)}$ and~$\matr{Y}_{21}\Emat^{(i)}$ are computed on-the-fly at each iteration~$i$ during the iterative solution of~\eqref{eq:finalsystem}, as follows.

Given~$\Emat^{(i)}$, the system of equations~\eqref{eq:AEFIEeqvdis} can be solved for ${\bigl[ jk_0\Hmat_{\mathrm{eq}}^{(i)} \quad c_0\rhomat_{\mathrm{eq}}^{(i)} \bigr]^T}$ at the~$i^{\mathrm{th}}$ iteration, which directly yields~$\matr{Y}_{11}\Emat^{(i)}$ and~$\matr{Y}_{21}\Emat^{(i)}$ via~\eqref{eq:invAEFIEeqvdis}.
For structures with multiple objects,~\eqref{eq:AEFIEeqvdis} is solved independently for each object.
For objects which are sufficiently small (less than $1\,000$ mesh triangles~\cite{TAPAIMin}), we use direct factorization.
For larger objects,~\eqref{eq:AEFIEeqvdis} is also solved with an iterative solver, and the AIM is applied independently for each object~\cite{TAPAIMin}, to compress the matrices~$\LmatA[l]$,~$\LmatPhi[l]$, and~$\Kpvmat[l]$.
In these cases, since the iterative solution of~\eqref{eq:AEFIEeqvdis} is ``nested'' into the iterative solution of~\eqref{eq:finalsystem}, we will refer to this as the ``nested'' system of equations.
An effective preconditioner is needed to ensure that the nested system~\eqref{eq:AEFIEeqvdis} is solved quickly at each iteration~$i$.
To this end, we use the constraint preconditioner proposed in~\cite{aefie2}.
If each object is relatively small compared to the full structure, the near-region entries of~$\LmatA[l]$ and~$\LmatPhi[l]$ can also be used in the preconditioner~\cite{near_zone_PC} to further reduce the iteration count for solving~\eqref{eq:AEFIEeqvdis}.
The results in~\secref{sec:results} confirm that the need to solve~\eqref{eq:AEFIEeqvdis} at each iteration~$i$ during the solution of~\eqref{eq:finalsystem} does not have a significant adverse effect on the overall simulation time.
In general, for each object, the number of ``nested'' iterations required for solving~\eqref{eq:AEFIEeqvdis} is expected to be comparable to the number of iterations that would be required, when the AEFIE~\cite{aefie2} is applied for a perfect conductor with an identical geometry.

The iterative solution of~\eqref{eq:finalsystem} also requires an effective preconditioner.
Due to the similarlity of~\eqref{eq:finalsystem} to the system of equations in the eAEFIE~\cite{eaefie01}, we use a sparse preconditioner similar to the one proposed in~\cite{eaefie01}, defined as
\begin{align}
	\matr{P} =
	{%
	\begin{bmatrix}
		\mathrm{diag}\,\LmatA[m] & \mathrm{diag}\,\Pxin & \matr{D}^T\,\mathrm{diag}\,\LmatPhi[l]\,\matr{B} \\
		\mu_r\, \mathrm{diag}\,\LmatA & \mathrm{diag}\,\Kmat & \epsilon_{c,r}^{-1}\,\matr{D}^T\,\mathrm{diag}\,\LmatPhi\,\matr{B} \\
		\matr{F}\matr{D} & \matr{0} & k_0^2\,\matr{I}
	\end{bmatrix}}\label{eq:pcext}
\end{align}
As described in~\cite{eaefie01}, the analytical block inverse of~$\matr{P}$ can be constructed using the Schur complement, so that the preconditioner can be applied on-the-fly during the~$i^{\mathrm{th}}$ iteration.
Consequently, only a small number of iterations is needed even for very challenging problems, and the convergence behaviour of the proposed formulation is nearly identical to that of the eAEFIE, as will be shown in~\secref{sec:results}.



%% file: results.tex
\section{Results}\label{sec:results}

The accuracy and performance of the proposed formulation are demonstrated here through comparisons with state-of-the-art techniques from the literature.
All simulations were performed single-threaded on a 3\,GHz Intel Xeon CPU with $256\,$GB of memory.
The GMRES iterative solver~\cite{gmres} available in the scientific computation library PETSc~\cite{petsc-web-page} was used to solve~\eqref{eq:AEFIEeqvdis} and~\eqref{eq:finalsystem}.
A relative residual norm of~$10^{-4}$ was used as the convergence tolerance for GMRES in all cases, for the solution of both~\eqref{eq:AEFIEeqvdis} and~\eqref{eq:finalsystem}.
The surface current distributions were plotted with scientific colormaps taken from~\cite{scicolourmaps}.

To expedite the computation of the MGF~\cite{MGF01} during the direct integration step, the quasistatic part of the MGF is extracted in the spectral domain, and added back in analytical form in the spatial domain~\cite{qse}.
This speeds up the computation of the Sommerfeld integrals~\cite{SI_PE} for the terms remaining after the quasistatic extraction.
To avoid computing these integrals for every source-test interaction, interpolation tables are precomputed.

\subsection{Dielectric Sphere in Free Space}\label{sec:results:sph}

\begin{figure}[t]
  \centering
  \includegraphics[width=\linewidth]{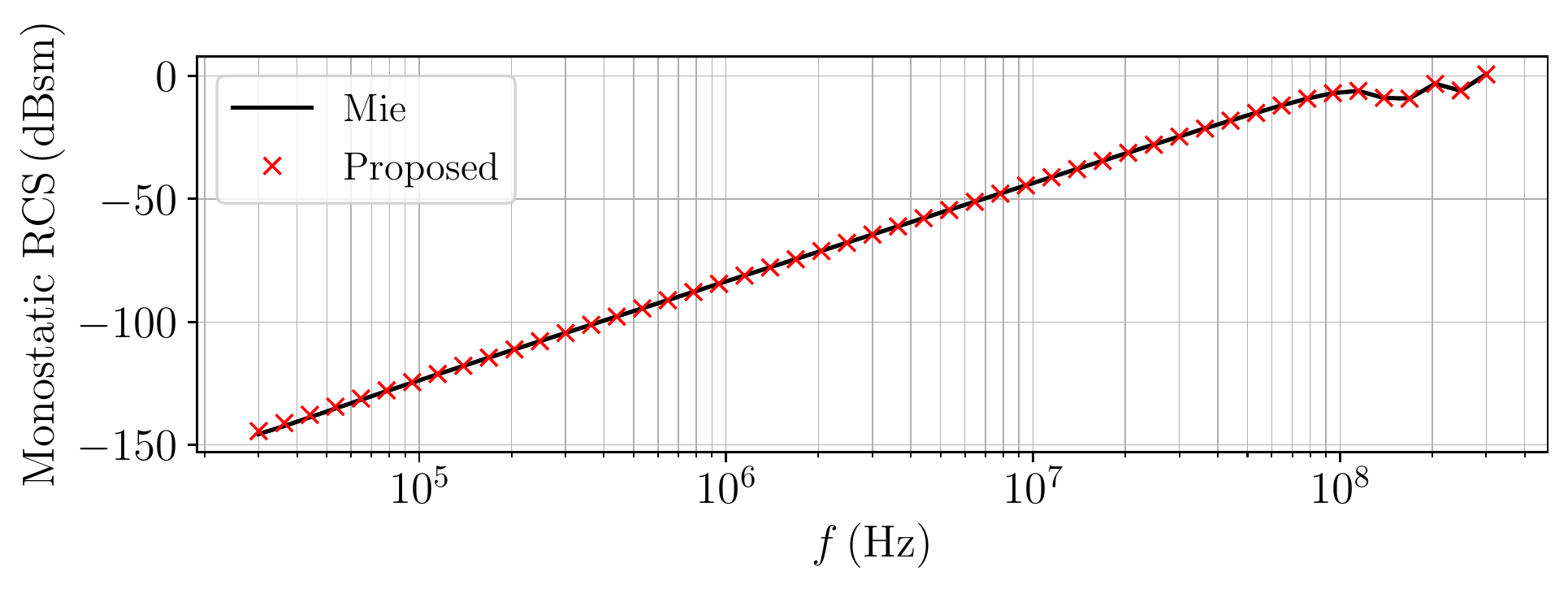}\\
  \includegraphics[width=\linewidth]{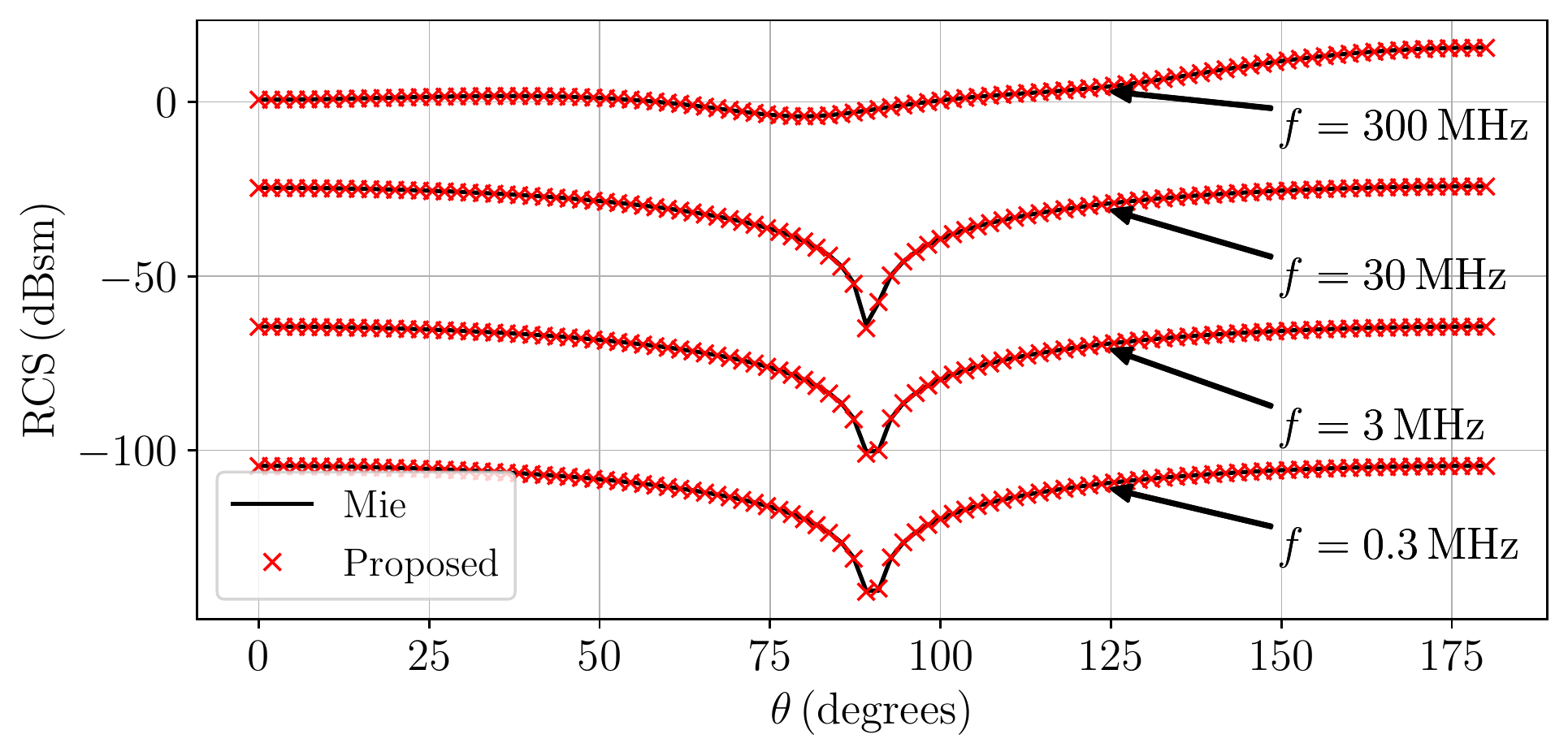}
  \caption{Accuracy validation for the sphere in \secref{sec:results:sph}. Top: monostatic RCS as a function of frequency. Bottom: bistatic RCS as a function of elevation angle, for~${\phi=0}$.}\label{fig:sph:acc}
\end{figure}

\begin{figure}[t]
	\centering
	\includegraphics[width=\linewidth]{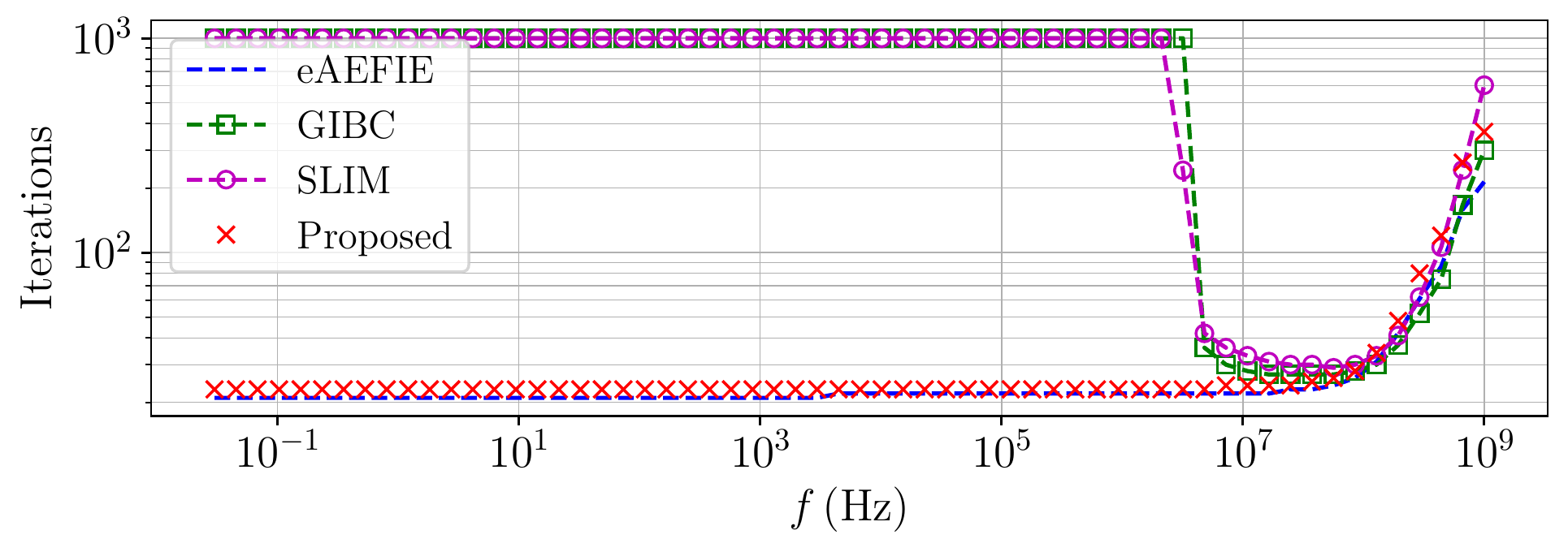}\\
	\includegraphics[width=\linewidth]{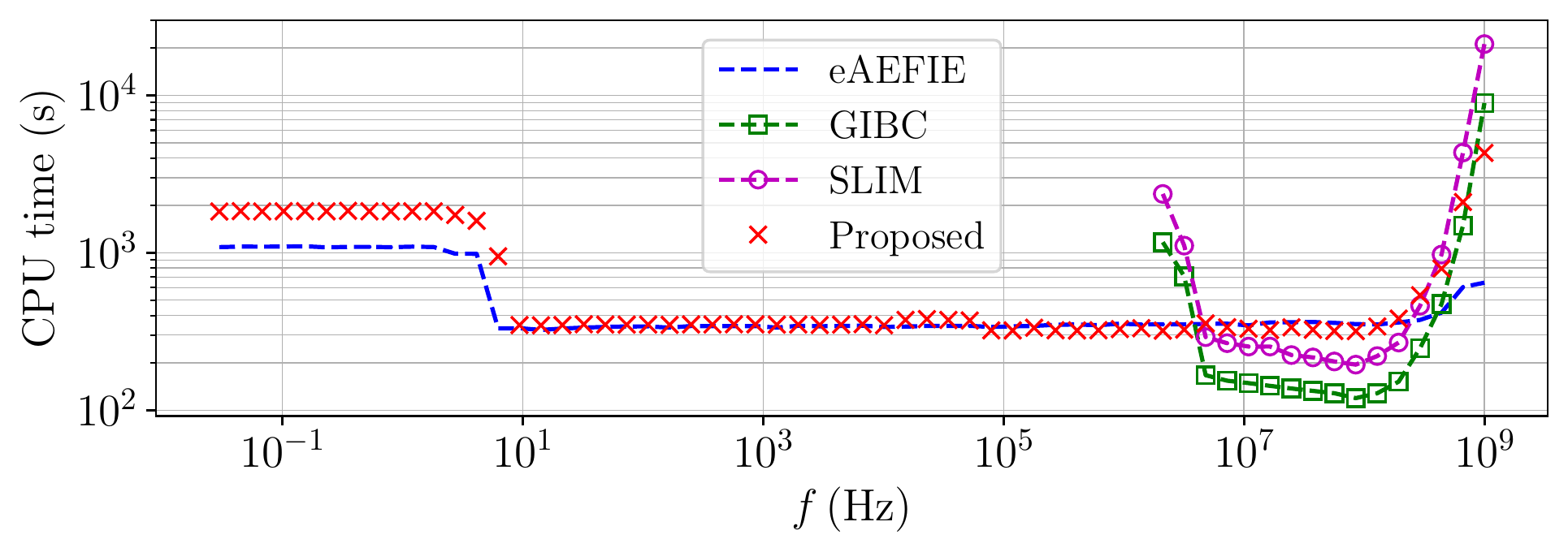}
	\caption{Performance comparison for the sphere in \secref{sec:results:sph}. Top: GMRES iterations. Bottom: total CPU time.}\label{fig:sph:perf}
\end{figure}

To validate the accuracy of the proposed formulation, we consider a dielectric sphere in free space, with relative permittivity~$12$ and diameter~$1\,$m.
The sphere is centered at the origin and meshed with~$3\,786$ triangles.
A plane wave, traveling in the~$-z$ direction with the electric field polarized along the~$x$ axis, is incident on the sphere.
The radar cross section (RCS) is measured and compared to the analytical solution obtained via the Mie series.
\figref{fig:sph:acc} demonstrates excellent agreement in both the monostatic and bistatic RCS between the proposed method and the Mie series.
The top panel of \figref{fig:sph:perf} shows the number of GMRES iterations required for convergence over a wide range of frequencies, for the proposed method compared to existing formulations.
Both the proposed and the eAEFIE~\cite{eaefie01} formulations show excellent convergence even at extremely low frequencies, while the GIBC~\cite{gibc} and SLIM~\cite{AWPLSLIM} formulations fail to converge within~$1\,000$ iterations, except above~${\sim}10^{5}\,$Hz.
The top panel of \figref{fig:sph:perf} also shows that the proposed method has the same limitation as the other AEFIE-based BEM techniques in that it becomes poorly conditioned at very high frequencies, as evidenced by the sharp increase in GMRES iterations.
This may be remedied with the use of more advanced preconditioners~\cite{HFEFIE,SIEPRECOND,RFCMP03}.
The bottom panel of \figref{fig:sph:perf} shows the total CPU time per frequency, which in this case is comparable to the CPU time for the eAEFIE, since the surrounding medium is homogeneous.
The fact that the GIBC and SLIM formulations avoid the use of BC basis functions is only advantageous over a relatively narrow frequency range, between approximately~$10^6\,$Hz and~$10^8\,$Hz.

\begin{table}[t]
	\captionsetup{width=0.93\linewidth}
	\centering
	\caption{Dielectric layer configurations for the numerical examples in~\secref{sec:results}. Layers are non-magnetic, and $h$ denotes the layer's height.}
	\begin{tabular}{ccc|cc|cc}
		\toprule
		\multicolumn{3}{c}{Sphere array} & \multicolumn{2}{c}{Inductor coil} & \multicolumn{2}{c}{SRR array} \\
		\multicolumn{3}{c}{(\secref{sec:results:spharr})} & \multicolumn{2}{c}{(\secref{sec:results:ind})} & \multicolumn{2}{c}{(\secref{sec:results:srr})} \\
		\cmidrule(lr){1-3}\cmidrule(lr){4-5}\cmidrule(lr){6-7}
		$\epsilon_r$ & $\sigma\,$(S/m) & $h$\,(mm) & $\epsilon_r$ & $h$\,($\mu$m) & $\epsilon_r$ & $h$\,($\mu$m) \\
		\midrule
		\multicolumn{2}{c}{Air} & \multicolumn{1}{c|}{$\infty$} & \multicolumn{1}{c}{Air} & \multicolumn{1}{c|}{$\infty$} & \multicolumn{1}{c}{Air} & \multicolumn{1}{c}{$\infty$}\\
		\midrule
		& & & $3.7$ & $50$ & $3.7$ & $5$ \\
		$3.47$ & $0.001725$ & $525$ & $11.9$ & $30$ & $4.4$ & $4$ \\
		& & & $4.4$ & $20$ & $2.1$ & $15$ \\
		\midrule
		\multicolumn{1}{c}{$6.5$} & \multicolumn{1}{c}{$0.003338$} & \multicolumn{1}{c|}{$\infty$} & \multicolumn{1}{c}{PEC} & \multicolumn{1}{c|}{$\infty$} & \multicolumn{1}{c}{PEC} & \multicolumn{1}{c}{$\infty$}\\
		\bottomrule
	\end{tabular}
	\label{tab:layers}
\end{table}

\subsection{Sphere Array Buried Underground}\label{sec:results:spharr}

\begin{figure}[t]
	\centering
	\includegraphics[width=\linewidth]{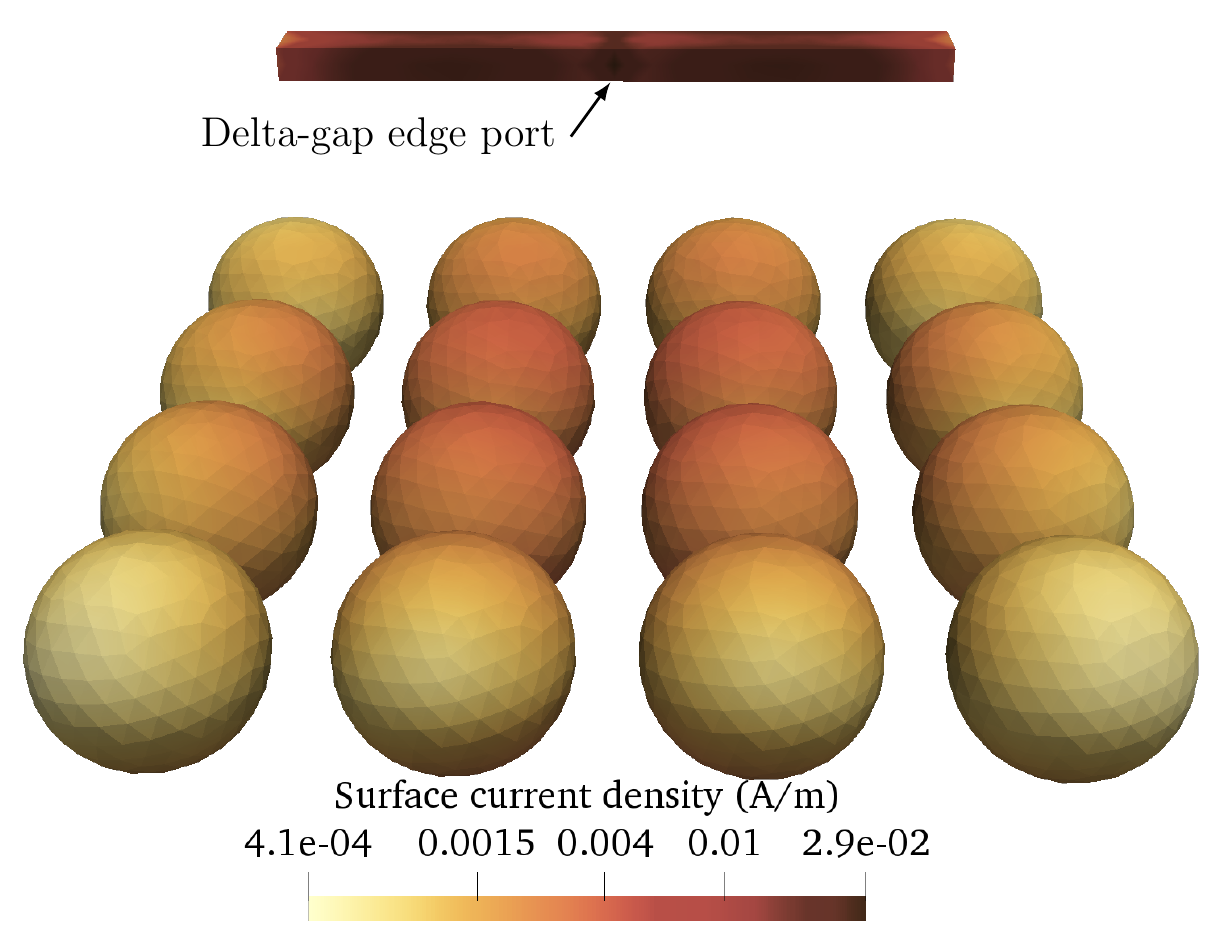}\\
	\caption{Geometry and electric surface current density at~$200\,$MHz for the array of spheres in \secref{sec:results:spharr}, excited by a dipole antenna.}\label{fig:spharrgeom}
\end{figure}

\begin{figure}[t]
	\centering
	\includegraphics[width=\linewidth]{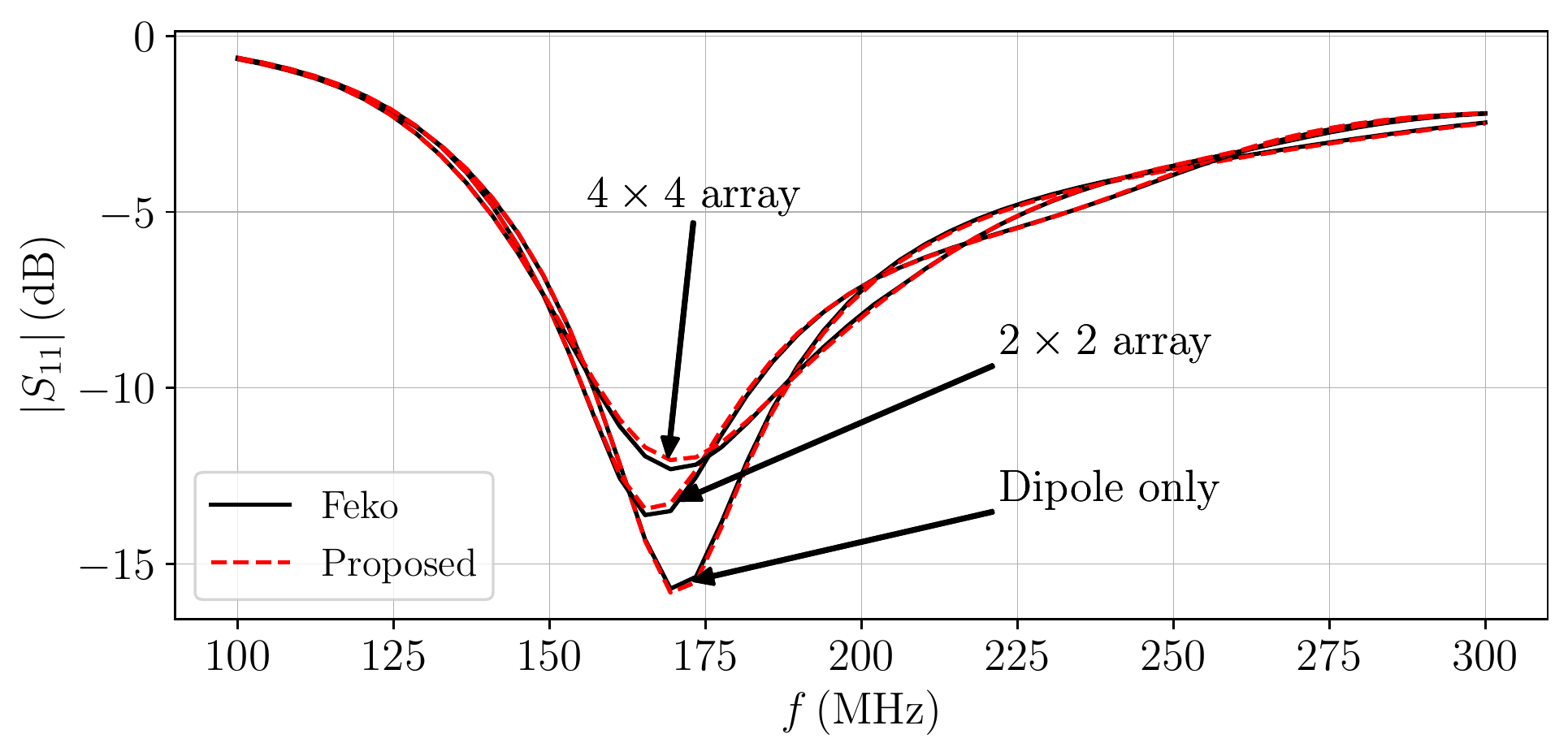}\\
	\includegraphics[width=\linewidth]{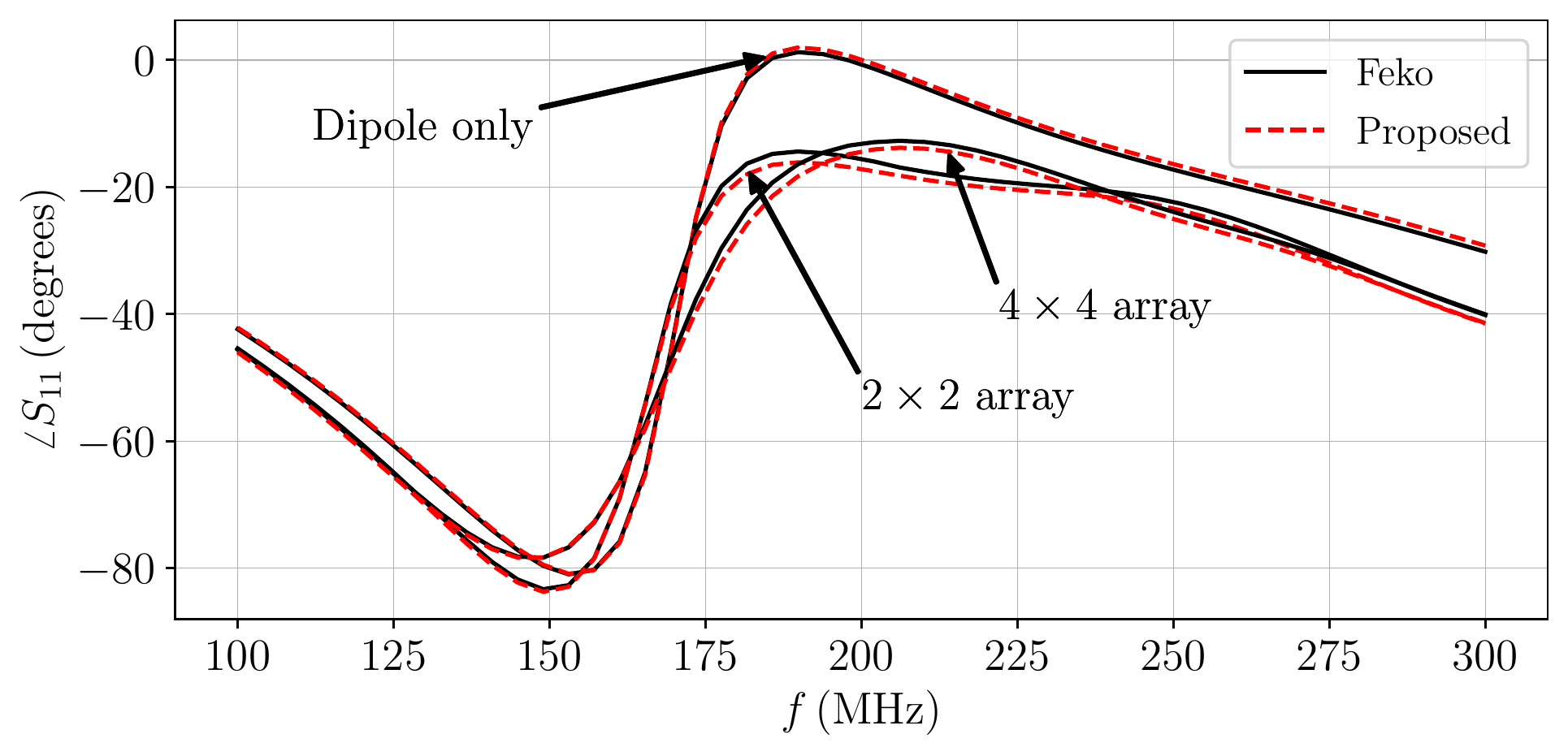}\\
	\caption{Scattering parameter validation for the array of spheres in \secref{sec:results:spharr}. Top: magnitude. Bottom: phase.}\label{fig:spharrS}
\end{figure}

\begin{figure}[t]
	\centering
	\includegraphics[width=\linewidth]{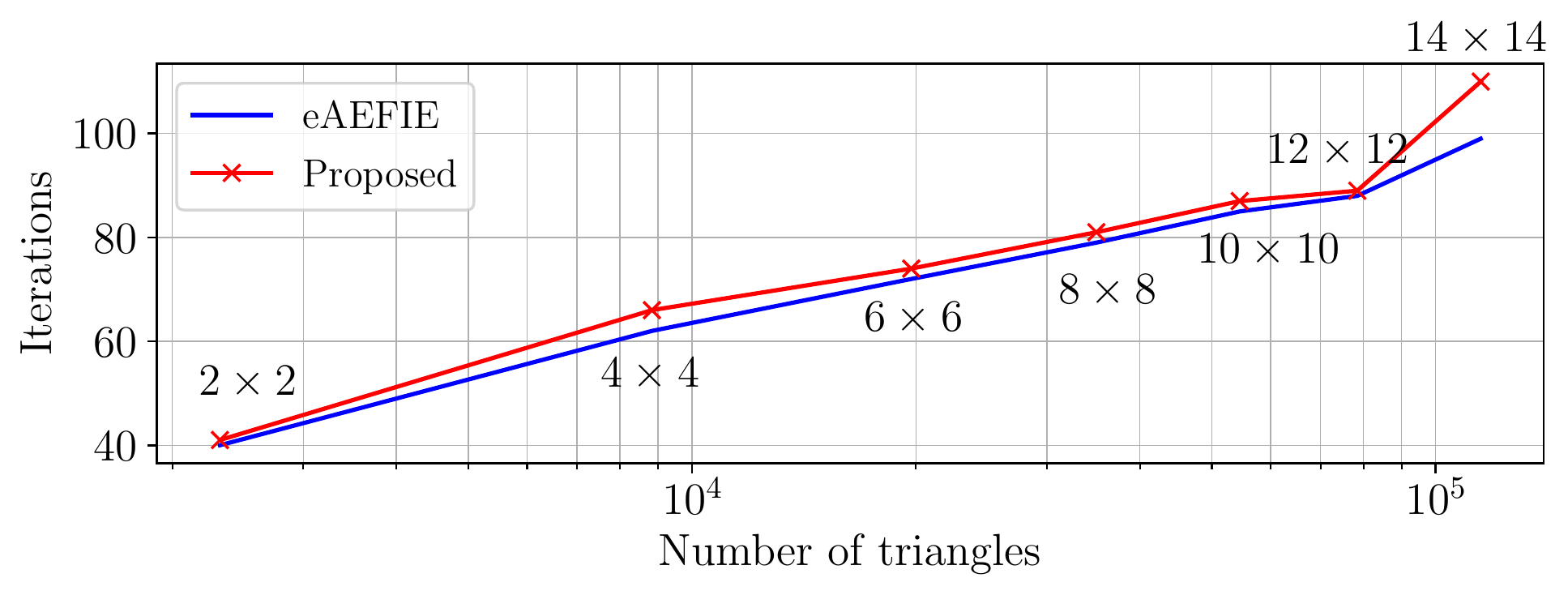}\\
	\includegraphics[width=\linewidth]{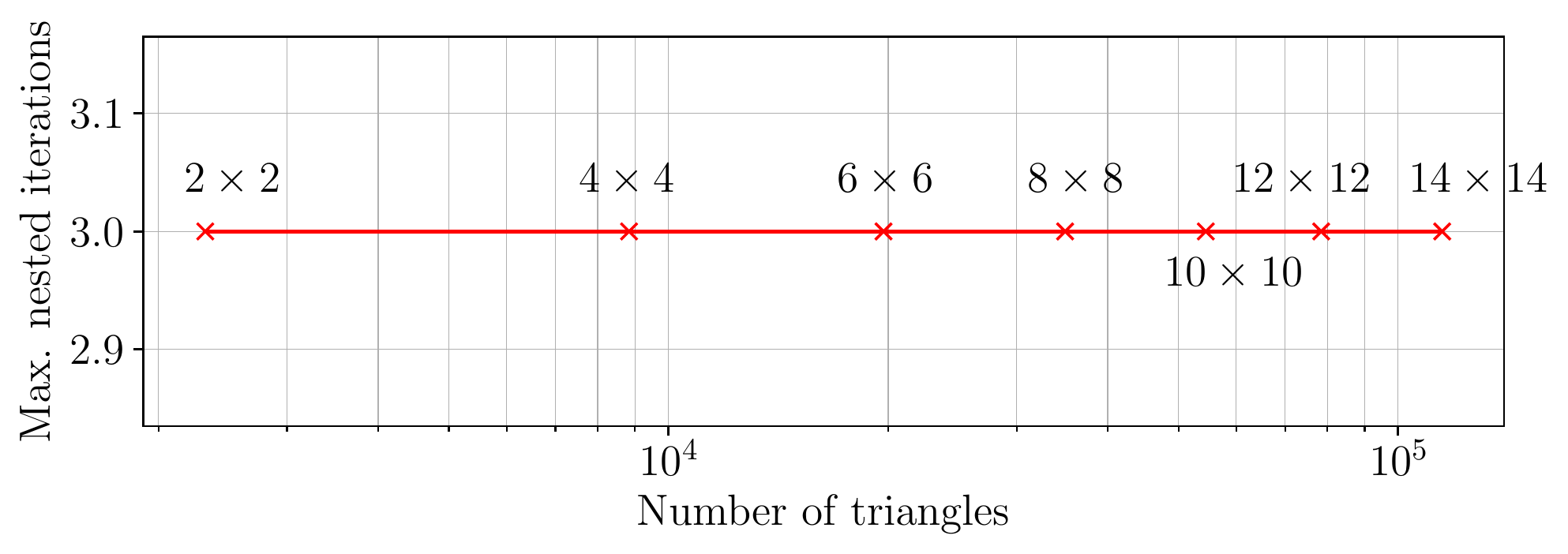}
	\caption{Number of GMRES iterations for the array of spheres in \secref{sec:results:spharr}. Top: solution of the final system~\eqref{eq:finalsystem}. Bottom: solution the ``nested'' system~\eqref{eq:AEFIEeqvdis} (worst case).}\label{fig:spharrNiter}
\end{figure}

\begin{figure}[t]
	\centering
	\includegraphics[width=\linewidth]{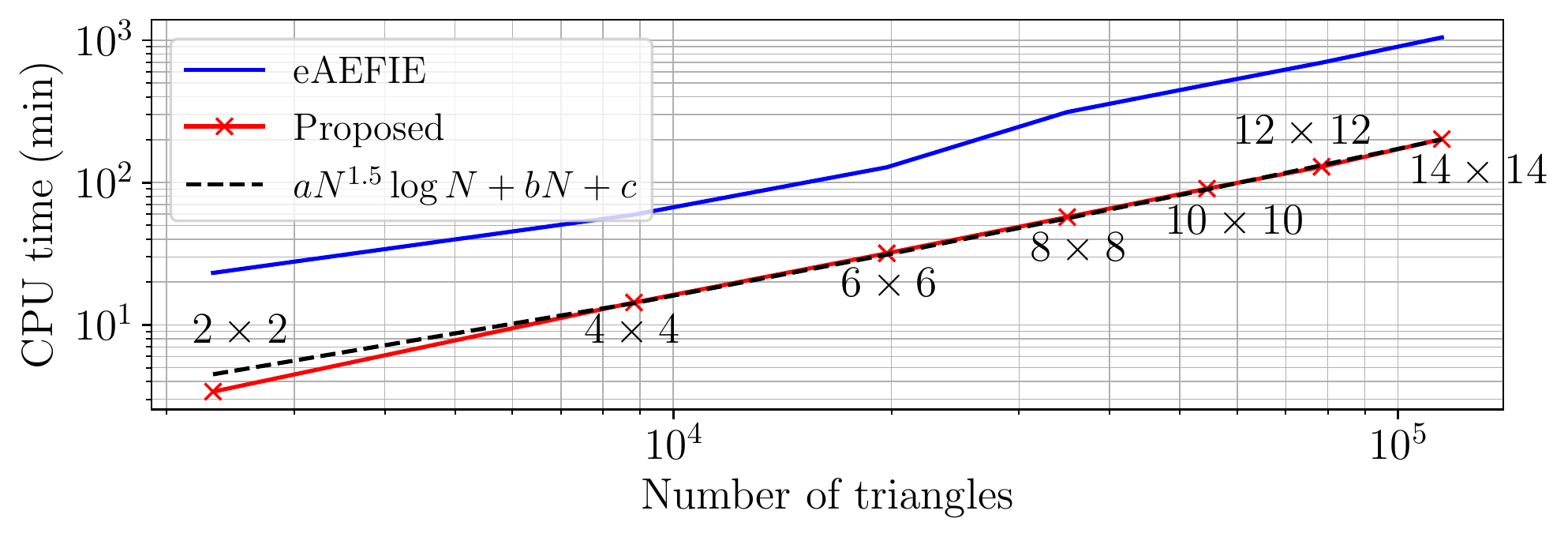}\\
	\includegraphics[width=\linewidth]{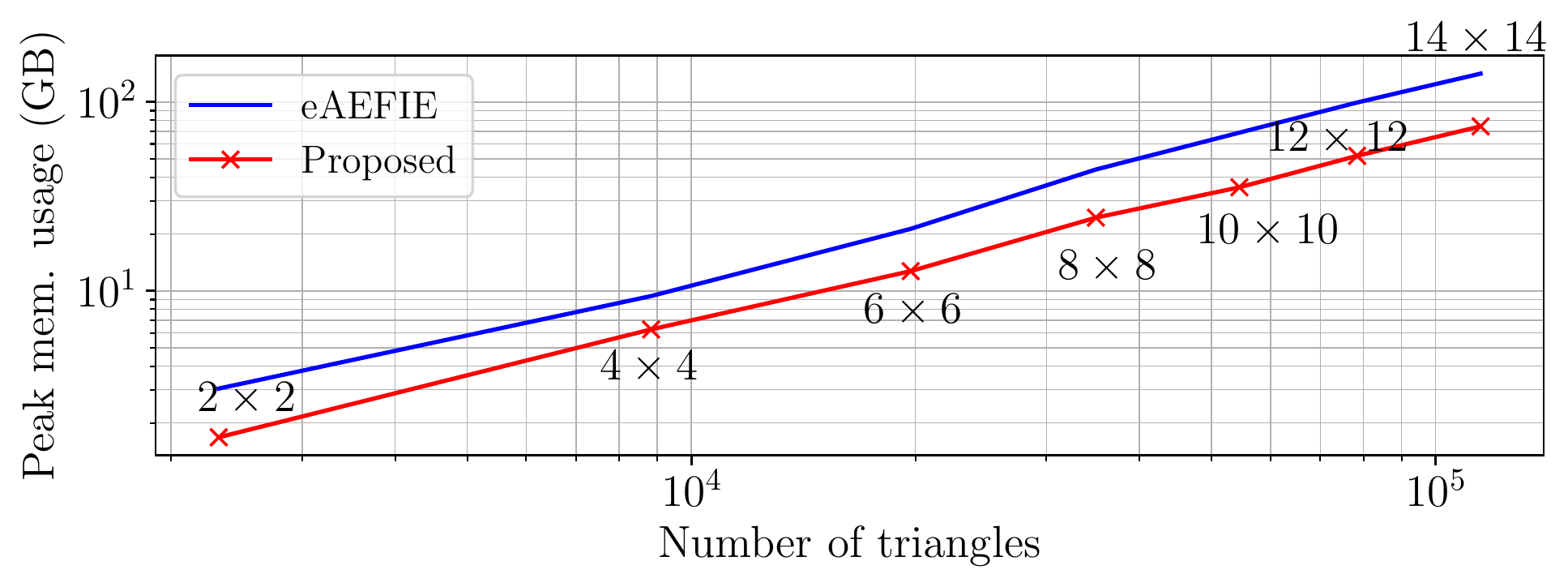}
	\caption{Performance comparison for the array of spheres in \secref{sec:results:spharr}. Top: total CPU time. Fit parameters are ${a=1.92 \times 10^{-7}}$, ${b=1.41 \times 10^{-3}}$, ${c=1.16}$, and $N$ is the number of triangles. Bottom: peak memory usage.}\label{fig:spharrprof}
\end{figure}

\begin{figure}[t]
  \centering
  \includegraphics[width=\linewidth]{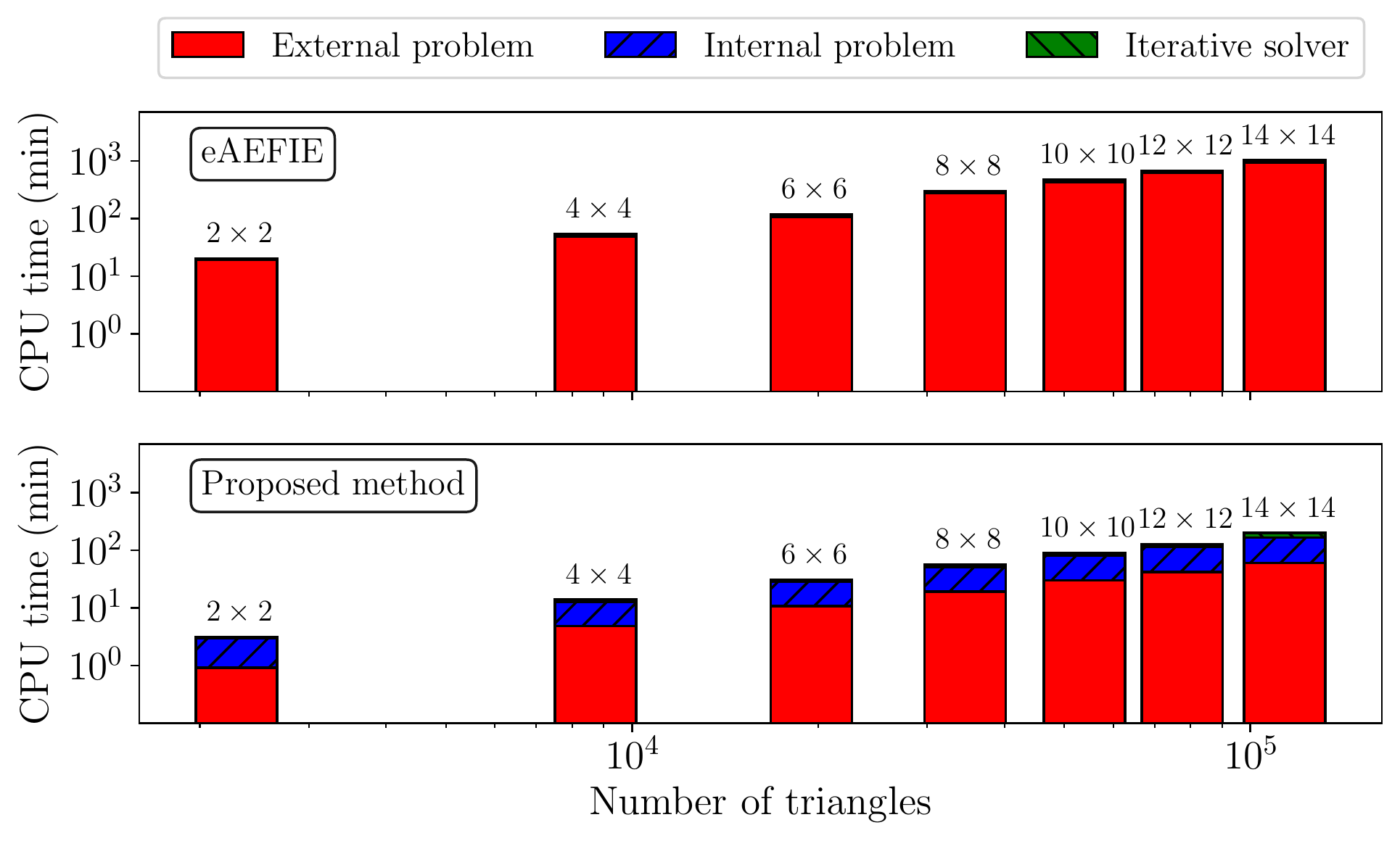}
  \caption{CPU time breakdown for the array of spheres in \secref{sec:results:spharr}. Top: eAEFIE. Bottom: proposed method.}\label{fig:spharrbrk}
\end{figure}

Next, we demonstrate the scalability of the proposed method in the context of a remote sensing application.
An array of dielectric spheres, each with relative permittivity~$12$ and diameter~$250\,$mm, is embedded in the layered medium described in the first column of \tabref{tab:layers}.
The considered dielectric layers represent typical electrical parameters of layered models of the ground~\cite{convcorr_AIM}.
We consider array sizes ranging from~${2 \times 2}$ to~${14\times 14}$ spheres, and the~${4 \times 4}$ case is visualized in \figref{fig:spharrgeom}.
In each case, the array is centered at the origin, and the bottom of the dielectric layer is located at~${z=-150}\,$mm.
A dipole antenna, designed to be half a wavelength long at~$200\,$MHz, is used to excite the array.
The antenna is modeled as a rectangular prism with dimensions~${75\,\text{mm} \times 750\,\text{mm} \times 37.5\,\text{mm}}$, and is centered at the point~$(0,0,493.75\,\text{mm})$.
A delta-gap edge port~\cite{gibson} is defined at the center of the dipole as shown in \figref{fig:spharrgeom}, where the scattering parameter~$S_{11}$ is measured.

\figref{fig:spharrS} shows the magnitude (top panel) and phase (bottom panel) of~$S_{11}$ in the vicinity of the antenna's design frequency, for the cases of~${2 \times 2}$ and~${4 \times 4}$ spheres. Also shown is the case with the antenna alone, in the absence of any spheres.
The results are compared against those obtained with the commercial solver Altair Feko.
Excellent agreement is observed over the entire bandwidth of the dipole in all cases.
The top panel of \figref{fig:spharrNiter} shows the GMRES iterations required for convergence as a function of the number of triangles for each sphere array, and demonstrates that the proposed method remains nearly as well conditioned as the eAEFIE, even as the problem size increases.
The bottom panel of \figref{fig:spharrNiter} shows the maximum number of ``nested'' GMRES iterations required for solving~\eqref{eq:AEFIEeqvdis}, and indicates the good conditioning of the system of equations associated with the internal problem in the equivalent configuration.
The top panel of \figref{fig:spharrprof} shows the total CPU time as a function of the number of triangles.
It demonstrates the significant improvement in performance of the proposed method compared to the eAEFIE, ranging from a~$6.8\times$ speed-up for the~${2 \times 2}$ array, to a~$5.2\times$ speed-up for the~${14 \times 14}$ array.
The fact that the double-layer operator is avoided in the external problem is primarily responsible for the efficiency of the proposed method compared to the eAEFIE.
Also shown in the top panel of \figref{fig:spharrprof} is a curve fit to the CPU time data for the proposed method, which confirms the predicted~$\mathcal{O}(N^{1.5}\log N)$ performance, where~$N$ is the number of mesh triangles.

The bottom panel of \figref{fig:spharrprof} shows the peak memory usage for the proposed method compared to the eAEFIE, and shows that in addition to CPU time savings, the proposed approach can also lead to a considerable reduction in memory usage, by up to a factor of two.
The reason for this reduction is that the global~$\Kpvmat[\mathrm{m}]$ operator of the eAEFIE, which contains entries for every pair of mesh edges in the entire structure, requires more memory than the \emph{relatively} sparse matrices~$\LmatA[l]$,~$\LmatPhi[l]$, and~$\Kpvmat[l]$ in the proposed method, which are local to each object.
The memory savings may be less pronounced when the structure contains large objects whose size is comparable to that of the entire structure.
In general, the memory usage of the proposed method is expected to be between~$0.5$ and~$1.5$ times that of the eAEFIE, depending on the structure.

Finally, \figref{fig:spharrbrk} shows, for the proposed method and the eAEFIE, a breakdown of the total CPU time into the three main contributions: the external problem, the internal problem, and the iterative solution of the final system of equations~\eqref{eq:finalsystem}.
The ``external problem'' bars include all near-region computations for matrix operators associated with the external region, which involve the MGF.
They also include the assembly and factorization cost associated with the preconditioner~\eqref{eq:pcext}.
The ``internal problem'' bars include all near-region computations for the object-wise matrix operators associated with the internal region, in both the original and equivalent configurations (where applicable).
For the proposed method, the ``internal problem'' bars also include the assembly and factorization cost associated with the constraint preconditioner used for solving~\eqref{eq:AEFIEeqvdis} iteratively, which is not applicable for the eAEFIE.
The ``iterative solver'' cost is the total time taken to solve~\eqref{eq:finalsystem} for a single excitation vector.
For the proposed method, this includes the cost of solving~\eqref{eq:AEFIEeqvdis} at each iteration.
\figref{fig:spharrbrk} clearly shows the dominant role played by computations related to the external region, and that the proposed method can significantly reduce this cost at the expense of a relatively small increase in the internal problem time.

\subsection{Inductor Coil with Microvias and a Dielectric Inclusion}\label{sec:results:ind}

\begin{figure}[t]
  \centering
  \includegraphics[width=\linewidth]{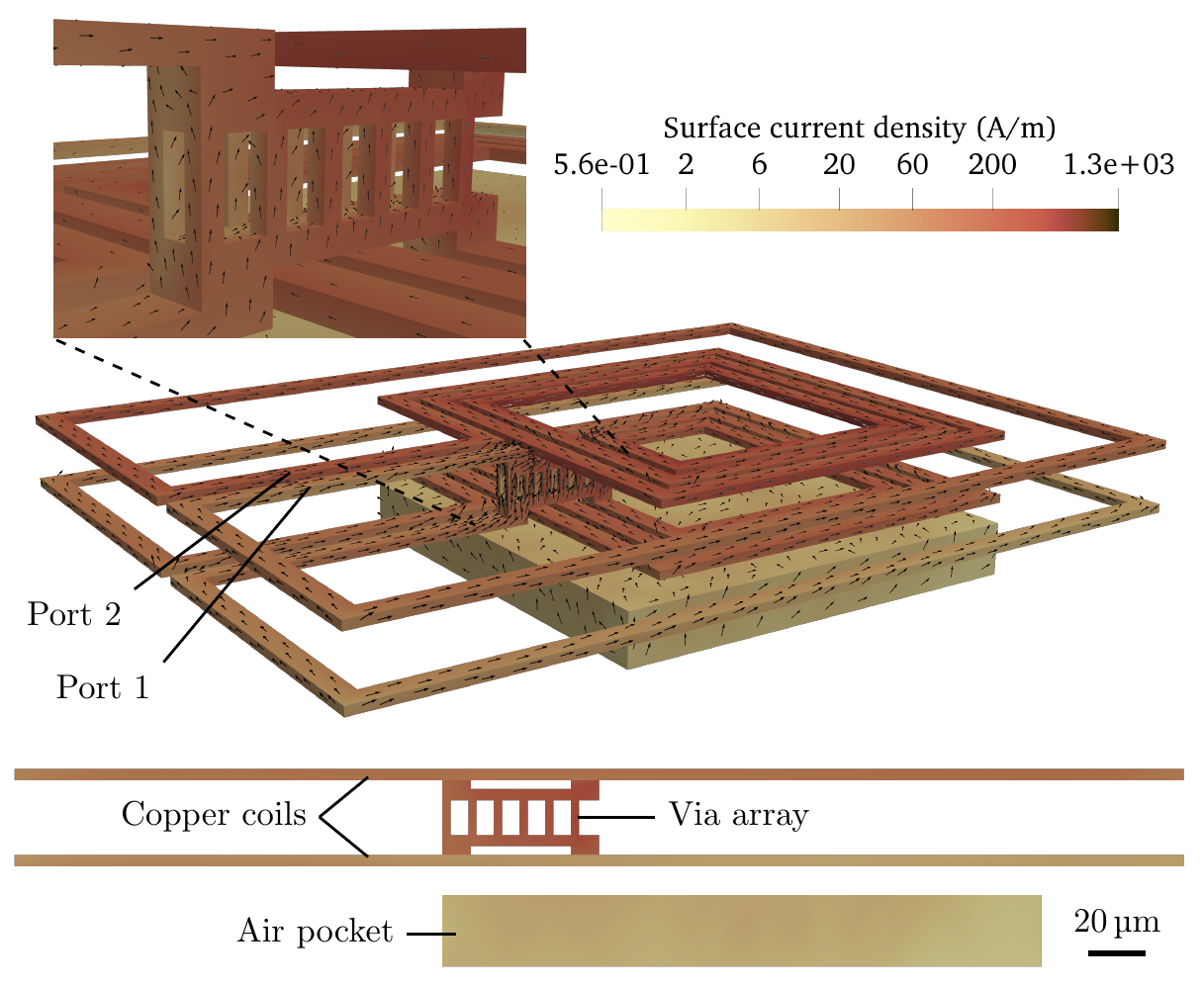}\\
  \caption{Geometry and electric surface current density at~$10\,$GHz for the inductor coil in \secref{sec:results:ind}. Inset shows a close-up of the microvia array.}\label{fig:indgeom}
\end{figure}

\begin{figure}[t]
  \centering
  \includegraphics[width=\linewidth]{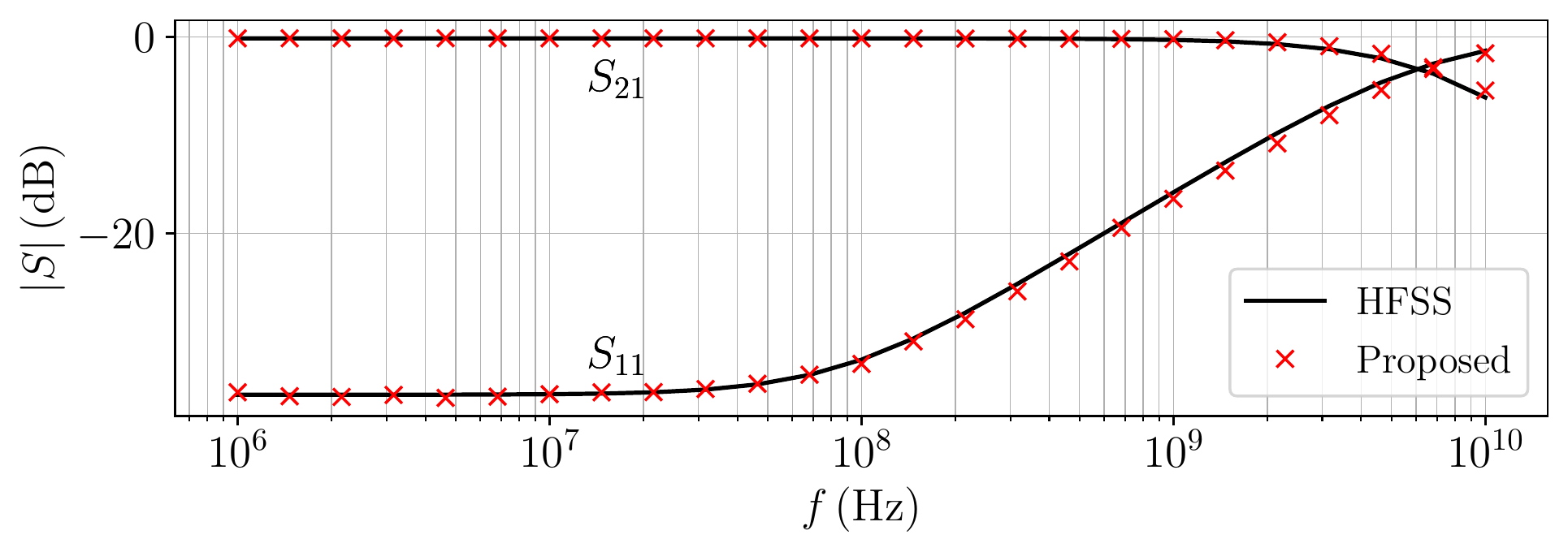}\\
  \includegraphics[width=\linewidth]{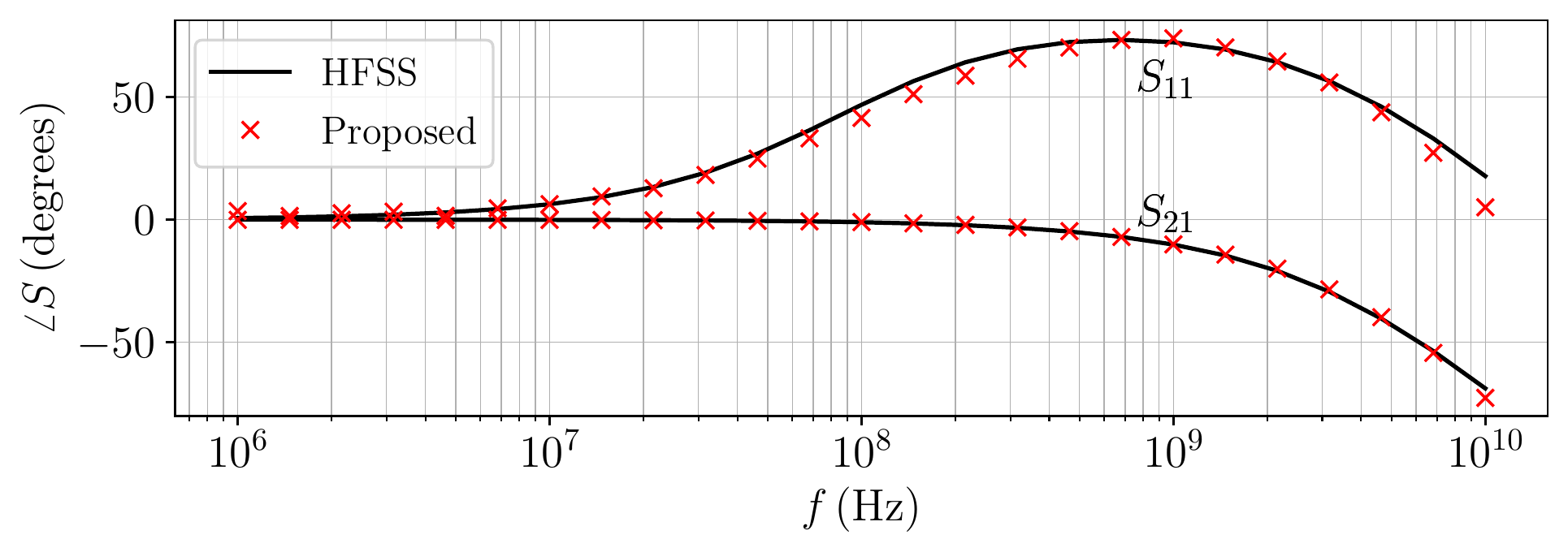}\\
  \caption{Scattering parameter validation for the inductor coil in \secref{sec:results:ind}. Top: magnitude. Bottom: phase.}\label{fig:indS}
\end{figure}

\begin{figure}[t]
  \centering
  \includegraphics[width=\linewidth]{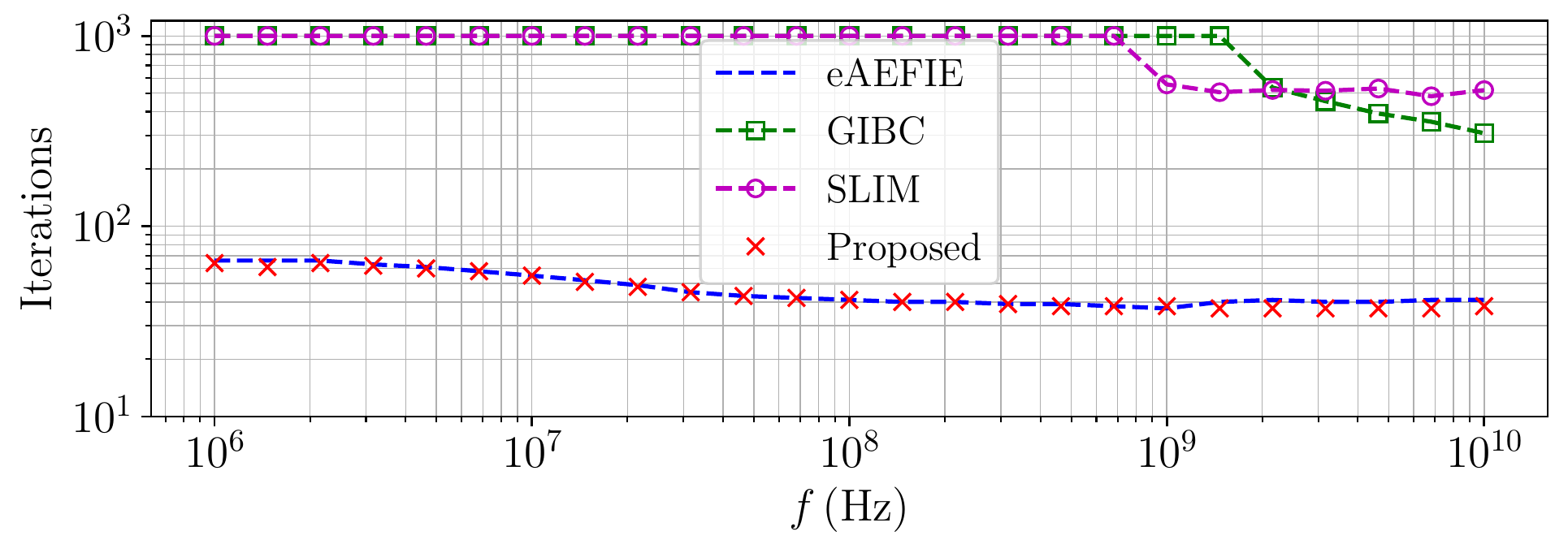}\\
  \includegraphics[width=\linewidth]{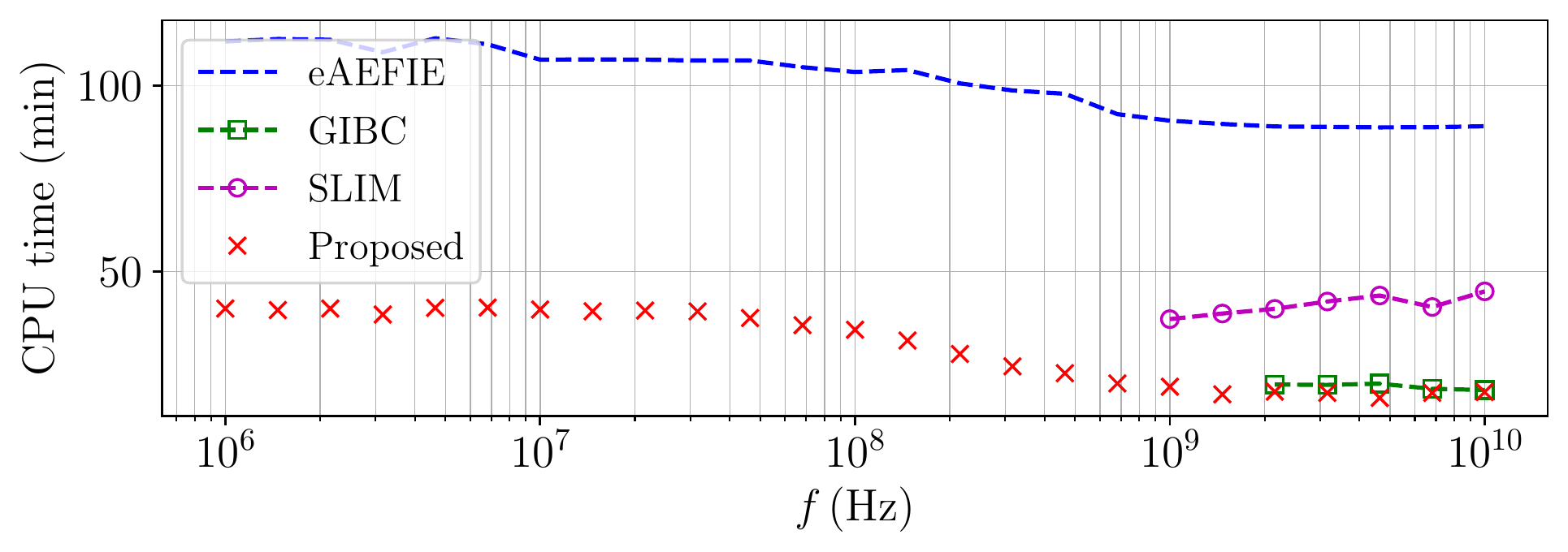}
  \caption{Performance comparison for the inductor coil in \secref{sec:results:ind}. Top: GMRES iterations. Bottom: total CPU time per frequency.}\label{fig:indprof}
\end{figure}

\begin{figure}[t]
  \centering
  \includegraphics[width=\linewidth]{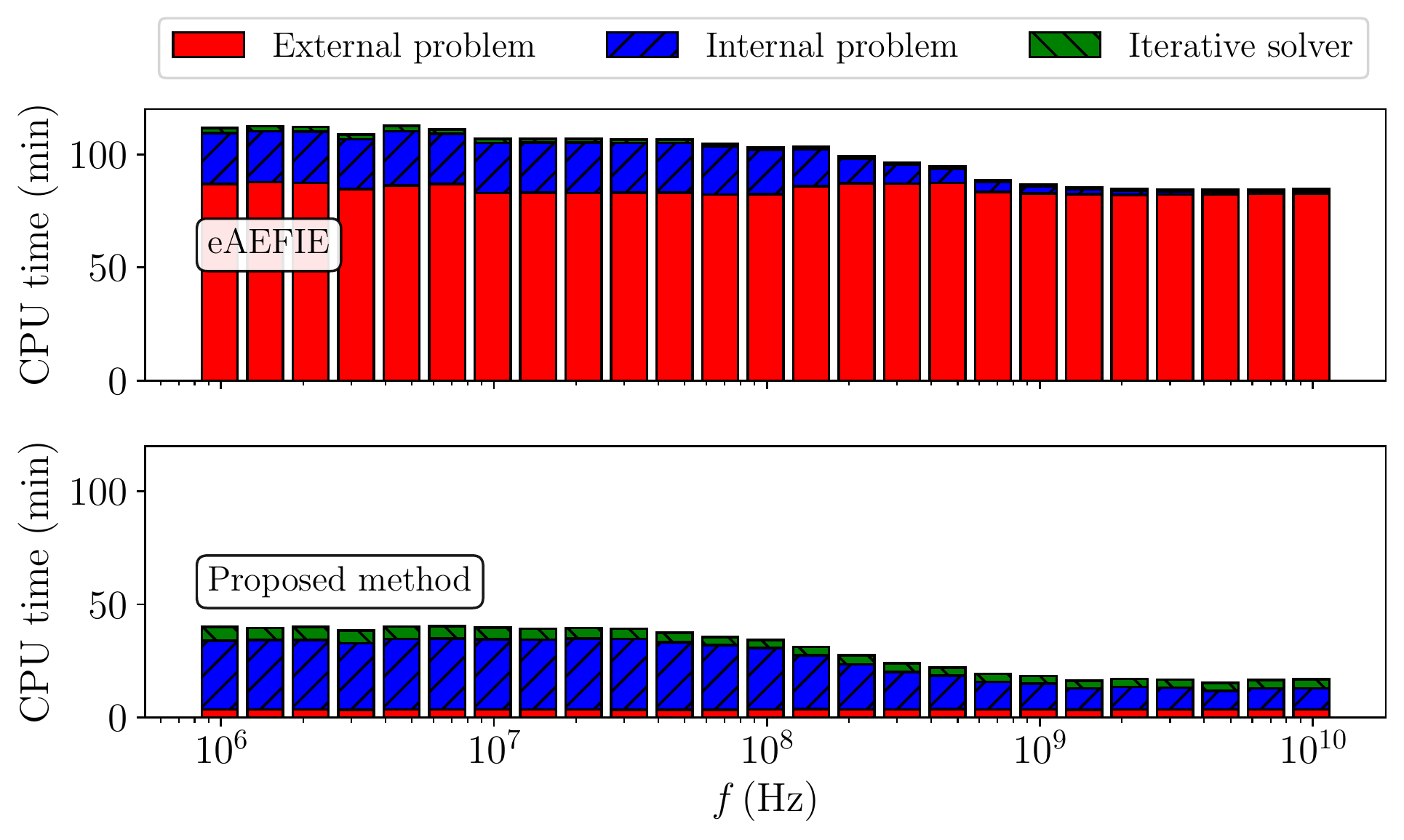}
  \caption{CPU time breakdown for the inductor coil in \secref{sec:results:ind}. Top: eAEFIE. Bottom: proposed method.}\label{fig:indbrk}
\end{figure}

Here, we apply the proposed formulation to analyze an on-chip inductor coil involving a microvia array, backed by an air pocket, shown in \figref{fig:indgeom}.
Air pockets are often introduced beneath inductor coils to reduce parasitic effects~\cite{indair01,indair02,indair03}.
Since air pockets are finite in size, they cannot be modeled as part of the background medium with the MGF.
Instead, the air pockets must be treated as dielectric objects whose surfaces are also meshed.
However, the proximity of a dielectric object to the metallic coils can lead to poor matrix conditioning in existing BEM formulations such as the GIBC~\cite{gibc} and SLIM~\cite{AWPLSLIM} approaches, as will be shown this section.
The presence of a microvia array further adds to the complexity of the problem, because the small features of the vias require a locally dense discretization.
Moreover, on-chip passive components must be characterized over a broad frequency range, over which the skin depth undergoes large variations.
Therefore, this is a multiscale problem in three respects: (a) material properties of adjacent objects, (b) relative electrical sizes of features within the structure at a particular frequency, and (c) variations in the skin depth and the structure's overall electrical size over the considered frequency range.

The structure consists of two copper inductor coils, each of which is a~$4\times$ scaled-up version of the coil described in~\cite{EPEPS2017}, which in turn is based on~\cite{fastmaxManual}.
The two coils are stacked back-to-back, and connected with a~${2 \times 6}$ array of vias.
A side view of the stacked coils is shown in the bottom panel of \figref{fig:indgeom}.
Each via has a rectangular cross section of size~\SI{3}{\micro\metre}\,$\times$\,\SI{3}{\micro\metre}, and a height of~\SI{12}{\micro\metre}.
The air pocket has dimensions of~\SI{210}{\micro\metre}\,$\times$\,\SI{210}{\micro\metre}\,$\times$\,\SI{25}{\micro\metre}, and is placed \SI{10}{\micro\metre} below the lower coil.
The background multilayer substrate is described in the second column of \tabref{tab:layers}.
The structure is meshed with $4\,758$ triangles, and is excited with two delta-gap edge ports~\cite{gibson}, as shown in \figref{fig:indgeom}.

\figref{fig:indS} shows the simulated scattering ($S$) parameters over a broad range of frequencies, from~$1\,$MHz to~$10\,$GHz.
The top panel shows the magnitude, while the bottom panel shows the phase.
The results are compared against those from a commercial finite element solver (Ansys HFSS 2020 R2).
Very good agreement is observed over the entire frequency range, and the development of skin effect is accurately captured.
The top panel of \figref{fig:indprof} shows the number of GMRES iterations required for convergence, for the proposed method compared with some representative existing BEM formulations.
The proposed method and the eAEFIE~\cite{eaefie02} show excellent convergence over the entire frequency range, while the GIBC~\cite{gibc} and SLIM~\cite{AWPLSLIM} formulations fail to converge for frequencies below~${\sim}1\,$GHz.
Even for frequencies above~${\sim}1\,$GHz, the GIBC and SLIM methods require approximately~$10\times$ more iterations.
Since each port requires invoking the iterative solver, the large number of iterations required by the GIBC and SLIM approaches would become a bottleneck for structures with tens or hundreds of ports.
At each frequency point, the proposed method required at most $35$ ``nested'' iterations for solving~\eqref{eq:AEFIEeqvdis}.

The bottom panel of \figref{fig:indprof} shows the total CPU time taken per frequency, and shows the excellent performance of the proposed formulation in comparison to existing methods.
The proposed method provides an overall~$3.4\times$ speed-up compared to the eAEFIE, by obviating the need for the double-layer operator in the external problem.
A maximum of~$7.7\,$GB of memory was used by the proposed approach over the entire sweep, while the eAEFIE used a maximum of~$9.4\,$GB.
\figref{fig:indbrk} shows, for the proposed method and the eAEFIE, a breakdown of the total CPU time similar to the one provided in \secref{sec:results:spharr}.
The drastic reduction in the external problem cost is clearly visible, and the computational bottleneck associated with the MGF is effectively eliminated, with a marginal increase the internal problem and iterative solver costs.
The advantage of the proposed method is more pronounced here than for the array of spheres in \secref{sec:results:spharr}, because the structure is denser and the number of elements in the near region is signficantly larger.
Therefore, in the eAEFIE, the near-region computations associated with~$\Kpvmat_{\mathrm{m}}$ dominate the CPU time, because~$\Kpvmat_{\mathrm{m}}$ involves both BC functions \emph{and} the MGF, while the proposed method uses no such operator.
\figref{fig:indbrk} indicates that even if there are multiple excitation vectors for which~\eqref{eq:finalsystem} must be solved iteratively, the proposed method can still provide a significant advantage over the eAEFIE.
The cost of the internal problem decreases with increasing frequency due to the development of skin effect, which sparsifies the internal problem matrices.

\subsection{Multiscale Split-Ring Resonator Array}\label{sec:results:srr}

\begin{figure}[t]
	\centering
	\includegraphics[width=\linewidth]{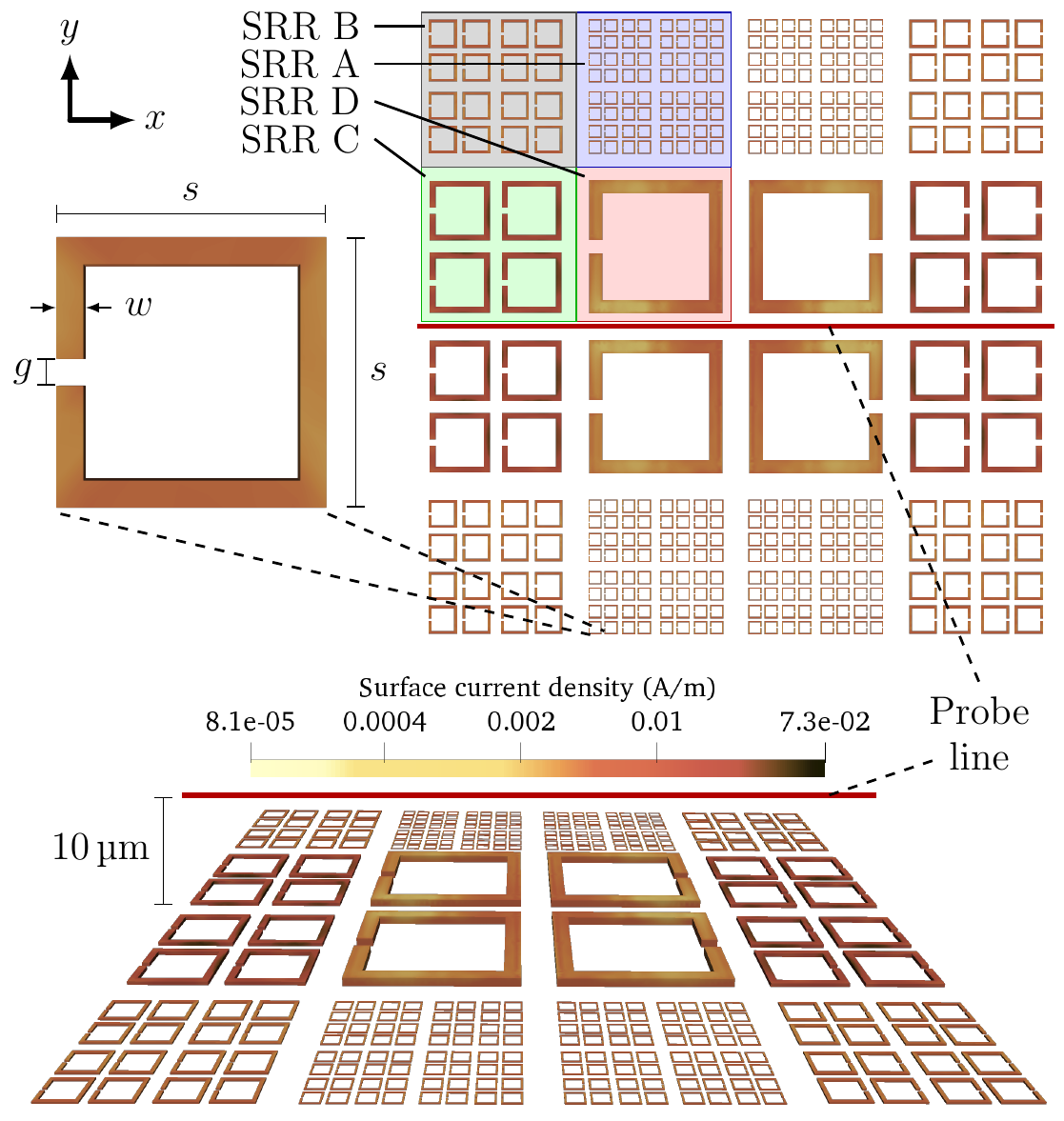}\\
	\caption{Geometry and electric surface current density at~$10\,$THz for the SRR array in \secref{sec:results:srr}.}\label{fig:srrgeom}
\end{figure}

\begin{table}[t]
	\captionsetup{width=0.96\linewidth}
	\centering
	\caption{Dimensions and electrical properties of each split ring in the SRR array in \secref{sec:results:srr}. Split rings are labeled based on the top panel of \figref{fig:srrgeom}, and $h$ denotes their height.}
	\begin{tabular}{ccccccc}
		\toprule
		& $s\,$(\SI{}{\micro\metre}) & $w\,$(\SI{}{\micro\metre}) & $g\,$(\SI{}{\micro\metre}) & $h\,$(\SI{}{\micro\metre}) & $\epsilon_r$ & $\sigma\,$(S/m) \\
		\midrule
		SRR A & $2$ & $0.2$ & $0.2$ & $0.1$ & $11$ & $0.4$ \\
		SRR B & $4.2$ & $0.42$ & $0.42$ & $0.21$ & $11$ & $0.4$ \\
		SRR C & $9$ & $0.9$ & $0.9$ & $0.45$ & $1$ & ${5.8 \times 10^7}$ \\
		SRR D & $20$ & $2$ & $2$ & $1$ & $6$ & $0$ \\
		\bottomrule
	\end{tabular}
	\label{tab:srrgeom}
\end{table}

\begin{figure}[t]
	\centering
	\includegraphics[width=\linewidth]{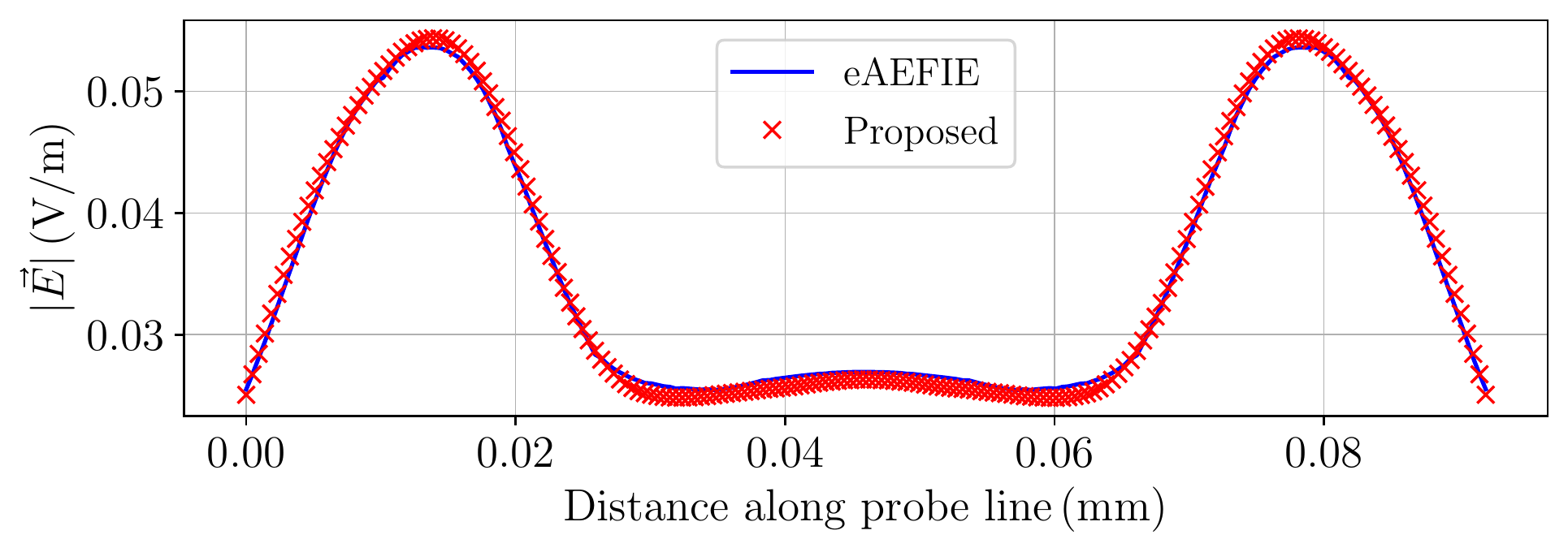}
	\caption{Magnitude of the electric field along a probe line for the SRR array in \secref{sec:results:srr}, at~$1\,$THz.}\label{fig:srrnf}
\end{figure}

\begin{figure}[t]
	\centering
	\includegraphics[width=\linewidth]{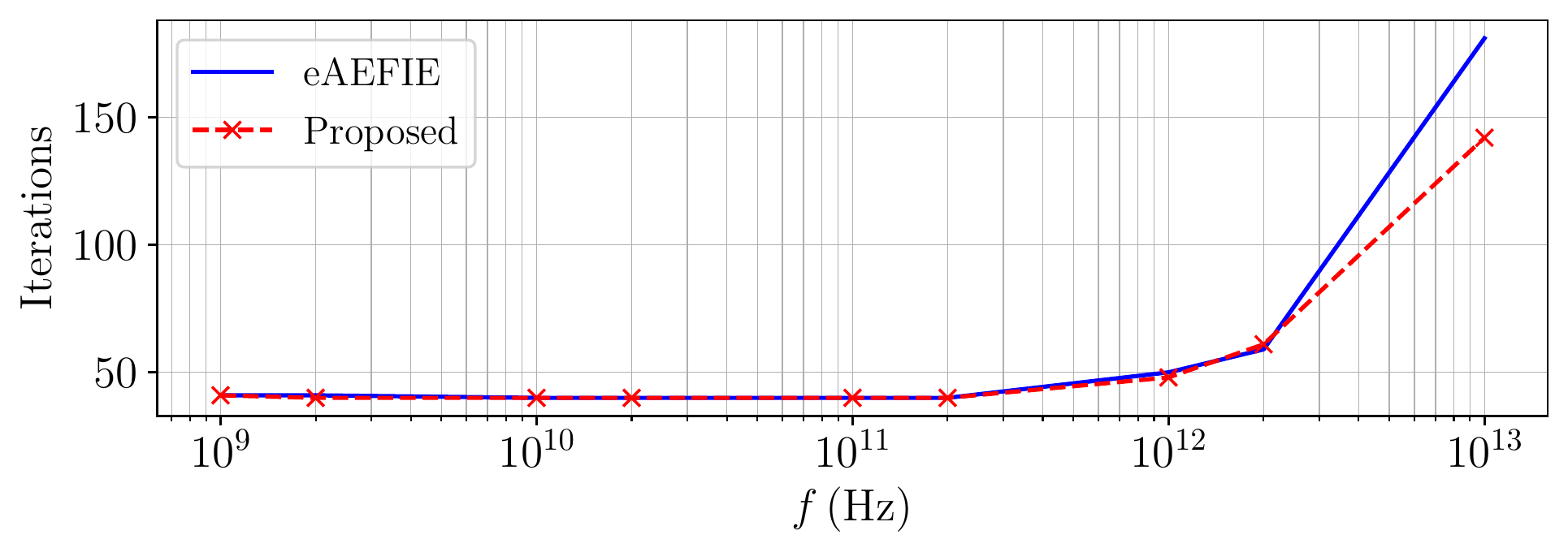}\\
	\includegraphics[width=\linewidth]{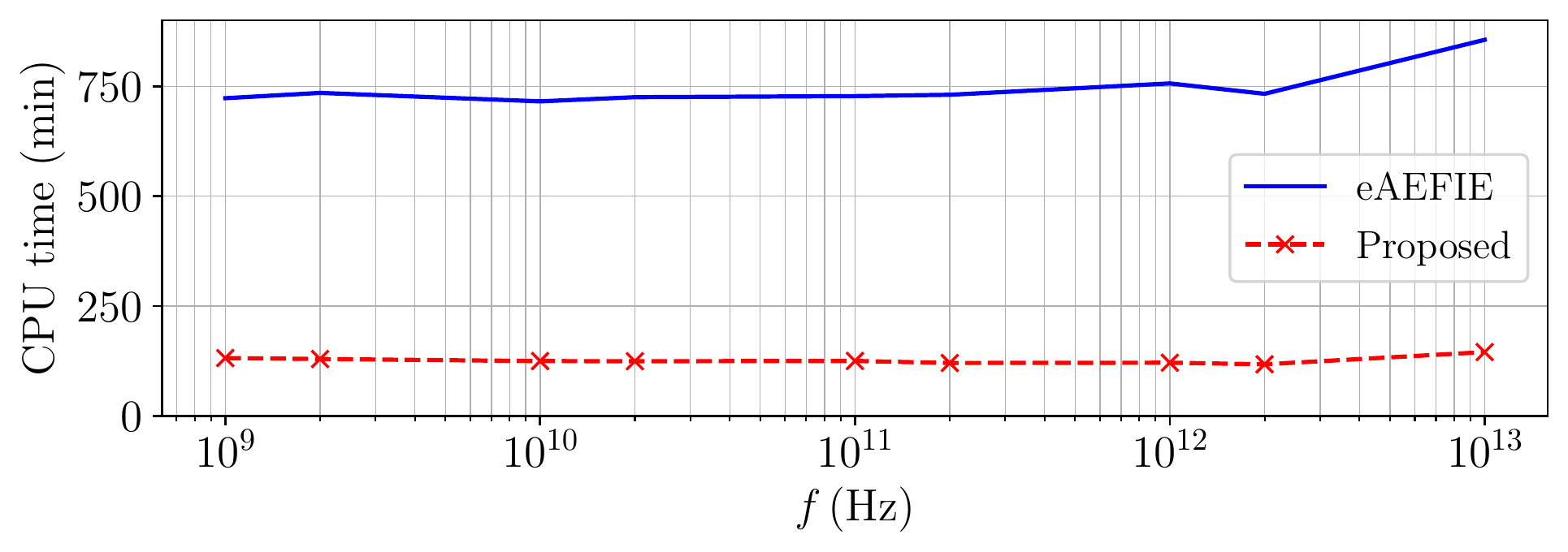}
	\caption{Performance comparison for the SRR array in \secref{sec:results:srr}. Top: GMRES iterations. Bottom: total CPU time per frequency.}\label{fig:srrprof}
\end{figure}

\begin{figure}[t]
  \centering
  \includegraphics[width=\linewidth]{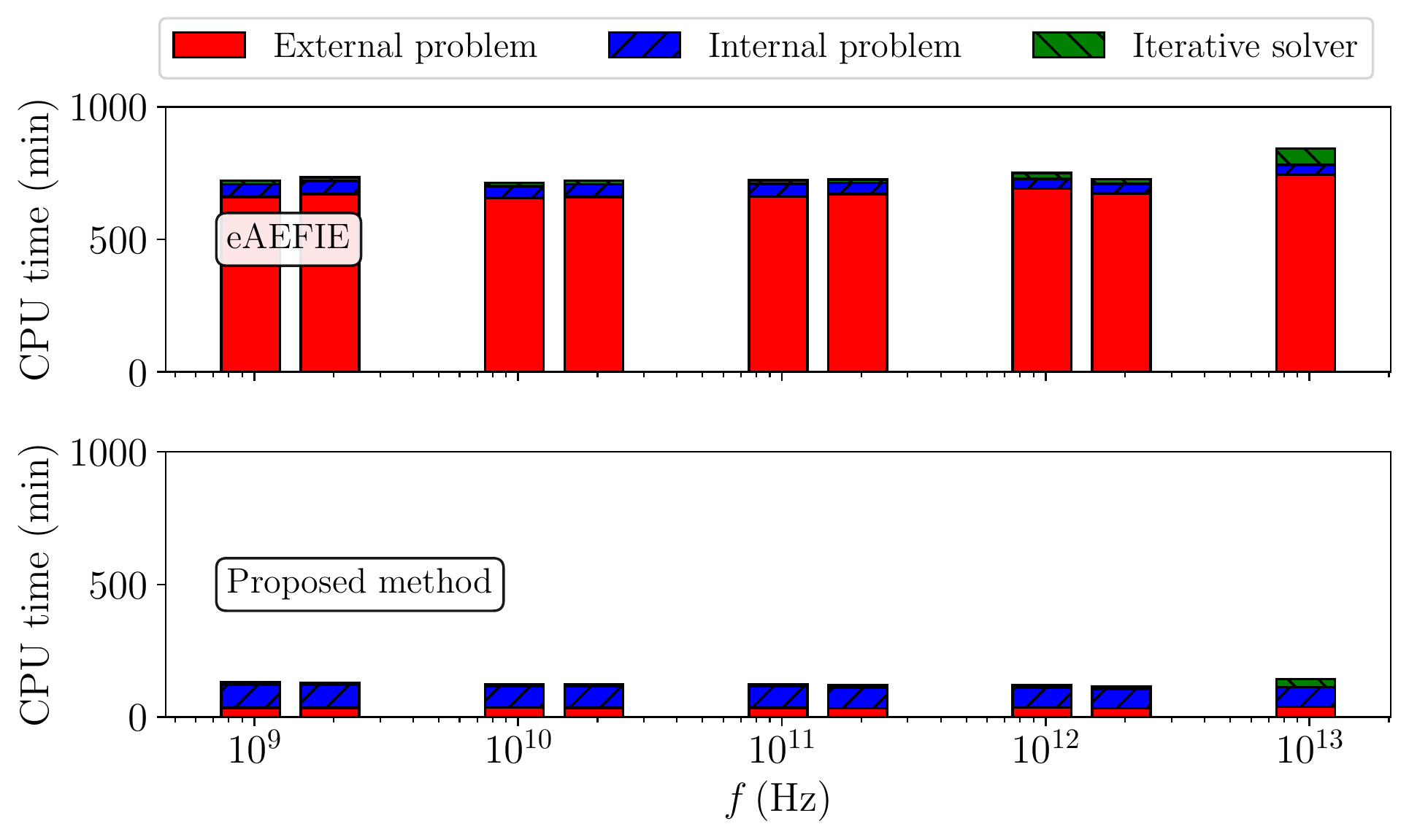}
  \caption{CPU time breakdown for the SRR array in \secref{sec:results:srr}. Top: eAEFIE. Bottom: proposed method.}\label{fig:srrbrk}
\end{figure}

Finally, we consider a multiscale array of dielectric and metallic split ring resonators (SRRs), shown in \figref{fig:srrgeom}.
SRRs are commonly used as unit cells in metasurface arrays, and can be made of dielectric or conductive materials~\cite{SRRmetaCond,optmeta01}.
The unit cells of metasurface arrays are typically sub-wavelength in size, while the entire array may span multiple wavelengths.
In addition, the arrays are typically fabricated on a substrate consisting of one or more dielectric layers.
Therefore, this example is meant to be representative of the computational challenges that may be encountered in the design and analysis of multiscale electromagnetic surfaces.

As shown in \figref{fig:srrgeom}, the structure contains four variations of unit cells.
The geometry of each unit cell is identical except for a scaling factor.
The top panel of \figref{fig:srrgeom} defines the pertinent dimensions of a single unit cell, and also shows the placement of each of the four variations of unit cells within the full array.
The variations are labeled, from the smallest to the largest unit cells, as ``SRR A'', ``SRR B'', ``SRR C'', and ``SRR D''.
\tabref{tab:srrgeom} provides the dimensions for each unit cell variation.
Some unit cells are made of dielectric materials, while others are conductive, and the material parameters are listed in \tabref{tab:srrgeom}.
The background layered medium is described in the third column of \tabref{tab:layers}.
The bottoms of all SRR unit cells are aligned at~${z=0}$, which coincides with the interface between the first (top-most) and second dielectric layers.
The structure is meshed with~$103\,568$ triangles, and is excited with a plane wave with the electric field oriented along the~$y$ axis, traveling in the~${-}z$ direction.

\figref{fig:srrnf} shows the electric field magnitude at~$1\,$THz for the proposed method, compared to the eAEFIE, measured along the probe line shown in \figref{fig:srrgeom}.
The proposed method is in good agreement with the eAEFIE.
The top panel of \figref{fig:srrprof} shows that the proposed method and the eAEFIE~\cite{eaefie01} converge within a reasonable number of GMRES iterations over a broad frequency range, from~$1\,$GHz to $10\,$THz.
The GIBC~\cite{gibc} and SLIM~\cite{AWPLSLIM} formulations did not converge within~$1\,000$ iterations for any of the frequency points considered.
We were unable to simulate this structure in HFSS or Feko within the available~$256\,$GB of memory.
Among all frequency points simulated, the proposed method required at most $23$ ``nested'' iterations for solving~\eqref{eq:AEFIEeqvdis}.
The bottom panel of \figref{fig:srrprof} shows the total CPU time per frequency, showing the significant computational advantage of the proposed formulation compared to the eAEFIE.
The proposed method yields an overall~$5.9\times$ speed-up compared to the eAEFIE formulation, reducing the total simulation time from~$4.6\,$days to~$18.9\,$hours.
The cost of computing the double-layer potential operator in the eAEFIE is particularly disadvantageous in this case, due to the intricate and dense nature of the structure.
In this case, the proposed method required $154\,$GB of memory at most, while the eAEFIE required $122\,$GB, because the structure contains several large unit cells.
The memory usage is still comparable between the two methods.
\figref{fig:srrbrk} shows the breakdown of the total CPU time per frequency, again demonstrating the significant impact of avoiding the double-layer operator for the external region.

In summary, the numerical tests considered in this section exemplify several applications where the proposed formulation can be a compelling alternative to existing techniques, such as the GIBC~\cite{gibc}, SLIM~\cite{AWPLSLIM} and eAEFIE~\cite{eaefie01,eaefie02} formulations.
The proposed method yields a well-conditioned system matrix while avoiding the double-layer potential operator in the external problem, unlike the eAEFIE, leading to a significant computational advantage for structures embedded in layered media.